\documentclass[reqno]{amsart}

\usepackage{amscd,amssymb}
\usepackage[cmtip,arrow,curve,matrix]{xy}

\usepackage{tikz-cd}

\newcommand{\splitpageyesno}[2]{#2}

\splitpageyesno{
\usepackage{splitpage}
}
{
\usepackage{fullpage}
}

\numberwithin{equation}{section}
\usepackage{setspace}

\usepackage{indentfirst}
\begin{document}
\title{On the inductive limit of direct sums of simple TAI algebras}
\author{Bo Cui, Chunlan Jiang, and Liangqing Li}

\maketitle

\makeatletter 
\renewcommand{\section}{\@startsection{section}{1}{0mm}
  {-\baselineskip}{0.5\baselineskip}{\bf\leftline}}
\makeatother

\section*{Abstract}

 \qquad~An ATAI (or ATAF, respectively) algebra, introduced in [Jiang1] (or in [Fa] respectively) is an inductive limit $\lim\limits_{n\rightarrow\infty}(A_{n}=\bigoplus\limits_{i=1}A_{n}^{i},\phi_{nm})$, where each $A_{n}^{i}$ is a simple separable nuclear TAI (or TAF) C*-algebra with UCT property. In [Jiang1], the second author classified all ATAI algebras by an invariant consisting orderd total K-theory and tracial state spaces of cut down algebras under an extra restriction that all element in $K_{1}(A)$ are torsion. In this paper, we remove this restriction, and obtained the classification for all ATAI algebras with the Hausdorffized algebraic $K_{1}$-group as an addition to the invariant used in [Jiang1]. The theorem is proved by reducing the class to the classification theorem of $\mathcal{AHD}$ algebras with ideal property which is done in [GJL1]. Our theorem generalizes the main theorem of [Fa] and [Jiang1] (see corollary 4.3).

\section{Introduction}
In [Lin2], Lin gave an abstract description of simple AH algebras (with no dimension growth) classified
in [EGL1]. He described the decomposition property of simple AH algebras in [G3] as TAI property and
proved that all simple separable nuclear TAI algebras with UCT are classifiable and therefore in the class of [EGL1]. For simple AH algebras of real rank zero, the corresponding decomposition property is called TAF
by [Lin1], which is partially inspired by Popa's paper [Popa].  As proved by Lin ([Lin1-2]), a simple separable nuclear C*-algebra $A$ with UCT property is a TAI (or TAF, respectively) algebra if and only if $A$ is a simple AH algebra (or simple AH algebra of real rank zero) with no dimension growth, which is classified in [EGL1].

As in [EG2], let $T_{II,k}$ be the 2-dimensional connected simplicial complex with $H_{1}(T_{II,k})=0$ and
$H_{2}(T_{II,k})=\mathbb{Z}/k\mathbb{Z}$, and let $I_{k}$ be the subalgebra of $M_{k}(C[0, 1])=C([0,1],M_{k}(\mathbb{C}))$ consisting of all functions $f$ with the properties $f(0)\in \mathbb{C}¡¤\mathbf{1}_{k}$ and $f(1)\in \mathbb{C}¡¤\mathbf{1}_{k}$ (this algebra is called an Elliott dimension drop interval algebra). Denote $\mathcal{HD}$ the class of algebras consisting of direct sums of building blocks of the forms $M_{l}(I_{k})$ and $PM_{n}(C(X))P$, with $X$ being one of the spaces $\{pt\}$, $[0, 1]$, $S^{1}$, and $T_{II,k}$, and with $P\in M_{n}(C(X))$ being a projection. (In [DG], this class is denoted by $SH(2)$, and in [Jiang1], this class is denoted by $\mathcal{B}$). A C*-algebra is called an $\mathcal{AHD}$ algebra, if it is an inductive limit of algebras in $\mathcal{HD}$. In [GJL2], the authors classified all $\mathcal{AHD}$ algebras with the ideal property.

In [Jiang1], Jiang classified ATAI algebras $A$ under the extra restriction that $K_{1}(A)=tor K_{1}(A)$. In this classification, Jiang used the scaled ordered total $K$-group (from [DG]) and the tracial state spaces $T(pAp)$ of cut-down algebras $pAp$, with certain compatibility conditions (from [Stev] and [Ji-Jiang]) as the invariant--we will call it $Inv^{0}(A)$. In this paper, we will use the invariant $Inv^{0}(A)$ together with the Hausdorffized algebraic $K_{1}$-group to deal with the torsion-free part of $K_{1}(A)$, with certain compatibility conditions--we will call this $Inv(A)$. We will prove that this invariant reduces to Jiang's invariant in the case that $K_{1}(A)$ is a torsion group, that is, we removed the extra restriction that $K_{1}(A)=tor K_{1}(A)$ of Jiang's classification in [Jiang1] and prove our theorem by reduced to the classification of $\mathcal{AHD}$ algebra with ideal property, which is done in [GJL1].

\section{Preliminaries And Definitions}

\subsection{2.1} Let $A$ be a C*-algebra. Two projections in $A$ are said to be equivalent if
they are Murray--von Neumann equivalent. We write $p \preceq q$ if $p$ is equivalent to a
projection in $qAq$. We denote by $[p]$ the equivalent class of projections equivalent
to $p$. Let $a \in A_{+}$, we write $p\preceq a$ if $p\preceq q$ for some projection $q \in aAa$.

\subsection{Definition 2.2} Let $G\subset A$ be a finite set and $\delta > 0$. We shall say that $\phi\in Map(A,B)$ is $G-\delta$ multiplicative if $\|\phi(ab)-\phi(a)\phi(b)\| < \delta$ for all $a, b \in G$.

We also use $Map_{G-\delta}(A,B)$ to denote all $G-\delta$ multiplicative maps.

If two maps $\varphi,\phi\in Map(A,B)$ satisfy the condition $\|\phi(a)-\psi(a)\|<\delta$ for all
$a\in G$, then we will write $\phi\approx_{\delta}\varphi$ on $G$.

If $G,H\subset A$ are two subsets of a C*-algebra $A$ and for any $g\in G$, there is a $b\in H$
with $\|g-b\|<\delta$, then we denote $G \subset_{\delta} H$.

\subsection{Definition 2.3}[Jiang1] We denote by $\mathcal{I}$ the class of all unital C*-algebras with the form $\bigoplus_{i=1}^{n}B_{i}$, where each $B_{i}\cong M_{k_{i}}$ or $B_{i}\cong M_{k_{i}}(C[0,1])$ for some $k_{(i)}$. Let $A\in \mathcal{I}$, we have the following:\\
(1) Every C*-algebra in $\mathcal{I}$ is of stable rank one.\\
(2) Two projections $p$ and $q$ in a C*-algebra $A\in \mathcal{I}$ are equivalent if and only if
$\tau(p)=\tau(q)$ for all $\tau\in T(A)$, where $T(A)$ denotes the space of all tracial states.\\
(3) For any $\varepsilon > 0$ and any finite subset $F\subset A$, there exist a number $\delta > 0$ and a finite subset $G \subset A$ satisfying the following: If $L : A \rightarrow B$ is a $G-\delta$ multiplicative contractive completely positive linear map, where $B$ is a C*-algebra, then there exists a homomorphism $h : A \rightarrow B$ such that
\[
\|h(a)-L(a)\|\leq \varepsilon, \quad \forall a\in F
\]

\subsection{Definition 2.4}([Lin1]) A unital simple C*-algebra in $A$ is said to be tracially AI (TAI) if
for any finite subset $F\subset A$ containing a nonzero element $b$,  $\varepsilon> 0$, integer $n > 0$
and any full element $a \in A_{+}$, there exist nonzero projection $p\in A$ and a C*-algebra
$I\subset A$ with $I\in \mathcal{I}$ and $1_{I}=p$ such that:\\
(1) $\|xp-px\|<\varepsilon$ for all $x \in F$;\\
(2) $pxp \in_{\varepsilon} I$ for all $x\in F$;\\
(3) $n[1-p]\preceq [p]$ and $1-p\preceq a$.

\subsection{Definition 2.5} [Jiang1] A C*-algebra $A$ (not necessary unital) is said to be ATAI algebra
(approximately TAI algebra) if it is the inductive limit of a sequence of direct sums
of simple unital TAI algebras with UCT.

\subsection{2.6} Let $A$ and $B$ be two C*-algebras. We use $Map(A,B)$ to denote the space of all completely positive $*$-contractions from $A$ to $B$. If both $A$ and $B$ are unital, then $Map(A,B)_{1}$ will denote the subset of $Map(A,B)$ consisting of all such unital maps.

\subsection{2.7} In the notation for an inductive limit system $\lim(A_{n},\phi_{n,m})$, we understand that
\[
\phi_{n,m}=\phi_{m-1,m}\circ\phi_{m-2,m-1}\circ\cdots\circ\phi_{n,n+1}
\]
where all $\phi_{n,m} : A_{n}\rightarrow A_{m}$ are homomorphisms.

We shall assume that, for any summand $A^{i}_{n}$ in the direct sum $A_{n}=\bigoplus_{i=1}^{t_{n}}A^{i}_{n}$, necessarily, $\phi_{n,n+1}(\mathbf{1}_{A^{i}_{n}})\neq 0$, since, otherwise, we could simply delete $A^{i}_{n}$ from $A_{n}$, without changing the limit algebra.

If $A_{n}=\bigoplus_{i}A^{i}_{n}$, $A_{m}=\bigoplus_{j}A^{j}_{m}$, we use $\phi^{i,j}_{n,m}$ to denote the partial map of $\phi_{n,m}$ from the i-th block $A^{i}_{n}$ of $A_{n}$ to the j-th block $A^{j}_{m}$ of $A_{m}$. Also, we use $\phi^{-,j}_{n,m}$ to denote the partial map of $\phi_{n,m}$ from $A_{n}$ to $A^{j}_{m}$. That is, $\phi^{-,j}_{n,m}=\bigoplus\limits_{i}\phi^{i,j}_{n,m}=\pi_{j}\phi_{n,m}$, where $\pi_{j}: A_{m}\rightarrow A^{j}_{m}$ is the canonical projection. We also use $\phi^{i,-}_{n,m}$ to denote $\phi_{n,m}|_{A^{i}_{n}}: A^{i}_{n}\rightarrow A_{m}$.

\subsection{2.8} An AH algebra is a C*-algebra which is  $\lim(A_{n}=\bigoplus^{k_{n}}_{i=1}p_{n,i}M_{[n,i]}(C(X_{n,i}))p_{n,i},\phi_{n,m})$, where each $X_{n,i}$ is a compact metrizable space, and $p_{n,i}\in M_{[n,i]}(C(X_{n,i}))$ is a projection. Recall in [G1], Gong proved that a simple AH algebra with uniformly bounded dimension of local spectra, i.e. $\sup_{n,i}\dim(X_{n,i})<\infty$ can be rewritten as AH inductive limit with spaces being $[0,1]$, $S^{1}$, $S^{2}$, $T_{II,k}$, $T_{III,k}$, where each $T_{II,k}$ (or $T_{III,k}$ respectively) is two--dimensional (or three--dimensional respectively) connected simplicial complexes with $H^{2}(T_{II,k})=\mathbb{Z}/k\mathbb{Z}$ and $H^{1}(T_{II,k})=0$ (or with $H^{3}(T_{III,k})=\mathbb{Z}/k\mathbb{Z}$ and $H^{1}(T_{III,k})=0$, $H^{2}(T_{III,k})=0$).

For any positive $k$ the dimension drop interval algebra $I_{k}$ is defined as
\[
I_{k}=\big\{f\in M_{k}(C[0,1])\big|f(0)=\lambda \mathbf{1}_{k},f(1)=\mu \mathbf{1}_{k}, \lambda,\mu\in \mathbb{C}\big\}.
\]

\subsection{2.9} [GJL1] Let $X$ be a compact space and $\psi: C(X)\rightarrow PM_{k_{1}}(C(Y))P$ $(rank(P)=k)$ be a unital homomorphism. For any point $y\in Y$, there are $k$ mutually orthogonal rank-1 projections $p_{1},p_{2},\cdots,p_{k}$ with $\sum\limits_{i=1}^{k}p_{i}=P(y)$ and $\big\{x_{1}(y),x_{2}(y),\cdots,x_{k}(y)\big\}\subset X$ (may be repeat) such that $\psi(f)(y)=\sum\limits_{i=1}^{k}f(x_{i}(y))p_{i}, \ \forall f\in C(X)$. We denote the set $\big\{x_{1}(y),x_{2}(y),\cdots,x_{k}(y)\big\}$ (counting multiplicities), by $Sp\psi_{y}$. We shall call $Sp\psi_{y}$ the spectrum of $\psi$ at the point $y$.

For any $f\in I_{k}$, let function $\underline{f}: [0,1]\rightarrow \mathbb{C}\bigsqcup M_{k}(\mathbb{C})$ (disjoint union) be defined by
\[
\underline{f}(t)=
\begin{cases}
\lambda,    & \     if \ \ t=0\  and  \ f(0)=\lambda\mathbf{1}_{k} \\
\mu,        & \     if \ \ t=1\  and  \ f(1)=\mu\mathbf{1}_{k} \\
f(t),       & \     if \ \ 0<t<1
\end{cases}
\]
That is, $\underline{f}(t)$ is the value of irreducible representation of $f$ corresponding to the point $t$. Similarly, for $f\in M_{l}(I_{k})$, we can define $\underline{f}: [0,1]\longrightarrow M_{l}(\mathbb{C}) \bigsqcup M_{lk}(\mathbb{C})$, by
\[
\underline{f}(t)=
\begin{cases}
a,    & \     if \ \ t=0\  and  \ f(0)=a\otimes \mathbf{1}_{k} \\
b,    & \     if \ \ t=1\  and  \ f(1)=b\otimes \mathbf{1}_{k} \\
f(t), & \     if \ \ 0<t<1
\end{cases}
\]

Suppose that $\phi: I_{k}\rightarrow PM_{n}(C(Y))P$ is a unital homomorphism. Let $r=rank(P)$. For each $y\in Y$, there are $t_{1},t_{2},\cdots,t_{m}\in [0,1]$ and a unitary $u\in M_{n}(\mathbb{C})$ such that
\[
P(y)=u\left(
        \begin{array}{cc}
          \mathbf{1}_{rank(P)} & 0 \\
                0     & 0 \\
        \end{array}
      \right)
u^{*}
\]
and
\[
(*) \qquad \qquad \ \ \phi(f)(y)=u\left(
  \begin{array}{ccccc}
    \underline{f}(t_{1})        &                 &               &                       &       \\
            &          \underline{f}(t_{2})       &               &                       &       \\
            &                   &              \ddots             &                       &       \\
            &                   &                 &          \underline{f}(t_{m})         &       \\
            &                   &                 &               &              \mathbf{0}_{n-r} \\
  \end{array}
\right)u^{*}\in P(y)M_{n}(\mathbb{C})P(y)\subset M_{n}(\mathbb{C})
\]
for all $f\in I_{k}$.

We define the set $Sp\phi_{y}$ to be the points $t_{1},t_{2},\cdots,t_{m}$ with possible fraction multiplicity. If $t_{i}=0$ or $1$, For example if we assume
\[
t_{1}=t_{2}=t_{3}=0<t_{4}\leq t_{5} \leq \cdots \leq t_{m-2}<1=t_{m-1}=t_{m},
\]
then $Sp\phi_{y}=\big\{0^{\sim\frac{1}{k}},0^{\sim\frac{1}{k}},0^{\sim\frac{1}{k}},t_{4},t_{5},\cdots,t_{m-2},1^{\sim\frac{1}{k}},1^{\sim\frac{1}{k}}\big\}$,
which can also be written as
\[
Sp\phi_{y}=\big\{0^{\sim\frac{3}{k}},t_{4},t_{5},\cdots,t_{m-2},1^{\sim\frac{2}{k}}\big\}
\]
Here we emphasize that, for $t\in (0,1)$, we do not allow the multiplicity of $t$ to be non-integral. Also for $0$ or $1$, the multiplicity must be multiple of $\frac{1}{k}$ (other fraction numbers are not allowed).

Let $\psi: C[0,1]\rightarrow PM_{n}(C(Y))P$ be defined by the following composition
\[
\psi: C[0,1]\hookrightarrow I_{k} \stackrel{\phi}{\longrightarrow} PM_{n}(C(Y))P,
\]
where the first map is the canonical inclusion. Then we have $Sp\psi_{y}=\{Sp\phi_{y}\}^{\sim k}$--that is, for each element $t\in (0,1)$, its multiplicity in $Sp\psi_{y}$ is exactly $k$ times of the multiplicity in $\phi_{y}$.

Recall that for $A=M_{l}(I_{k})$, every point $t\in (0,1)$ corresponds to an irreducible representation $\pi_{t}$, defined by $\pi_{t}(f)=f(t)$. The representations $\pi_{0}$ and $\pi_{1}$ defined by
\[
\pi_{0}=f(0) \qquad \ and \ \qquad  \pi_{1}=f(1)
\]
are no longer irreducible. We use $\underline{0}$ and $\underline{1}$ to denote the corresponding points for the irreducible representations. That is,
\[
\pi_{\underline{0}}(f)=\underline{f}(0) \qquad \ and \ \qquad  \pi_{\underline{1}}(f)=\underline{f}(1).
\]
Or we can also write $\underline{f}(0)\triangleq f(\underline{0})$ and $\underline{f}(1)\triangleq f(\underline{1})$. Then the equation $(*)$ could be written as
\[
\phi(f)(y)=u\left(
  \begin{array}{ccccc}
         f(t_{1})        &           &             &             &          \\
            &         f(t_{2})       &             &             &          \\
            &            &         \ddots          &             &          \\
            &            &           &          f(t_{m})         &          \\
            &            &           &             &       \mathbf{0}_{n-r} \\
  \end{array}
\right)u^{*}
\]
where some of $t_{i}$ may be $\underline{0}$ or $\underline{1}$. In this notation, $f(0)$ is equal to diag$\big( \underbrace{f(\underline{0}),f(\underline{0}),\cdots,f(\underline{0})}_{k}\big)$ up to unitary
equivalence.

Under this notation, we can also write $0^{\sim \frac{1}{k}}$ as $\underline{0}$. Then the example of $Sp\phi_{y}$ can be written as
\[
Sp\phi_{y}=\big\{0^{\sim \frac{1}{k}},0^{\sim \frac{1}{k}},0^{\sim \frac{1}{k}},t_{4},t_{5},\cdots,t_{m-2},1^{\sim\frac{1}{k}},1^{\sim\frac{1}{k}}\big\}=\big\{\underline{0},\underline{0},\underline{0},t_{4},t_{5},\cdots,t_{m-2},\underline{1},\underline{1}\big\}
\]

\subsection{2.10} We use $\mathcal{HD}$ to denote all C*-algebras $C=\bigoplus C^{i}$, where each $C^{i}$ is
of the forms $M_{l}(I_{k})$ or $PM_{n}C(X)P$ with $X$ being one of the spaces $\{pt\}$, $[0,1]$, $S^{1}$, $T_{II,k}$. Each block $C^{i}$ will be called a basic $\mathcal{HD}$ block or a basic building block.

By $\mathcal{AHD}$ algebra, we mean the inductive limit of
\[
A_{1}\stackrel{\phi_{1,2}}{\longrightarrow} A_{2} \stackrel{\phi_{2,3}}{\longrightarrow} \cdots \stackrel{}{\longrightarrow} A_{n}\cdots \stackrel{}{\longrightarrow}\cdots,
\]
where $A_{n}\in \mathcal{HD}$ for each $n$.

For an $\mathcal{AHD}$ inductive limit $A=\lim(A_{n},\phi_{n,m})$, we write $A_{n}=\bigoplus^{t_{n}}_{i=1}A^{i}_{n}$, where $A^{i}_{n}=P_{n,i}M_{[n,i]}(C(X_{n,i}))P_{n,i}$ of $A^{i}_{n}=M_{[n,i]}(I_{k_{n,i}})$. For convenience, even for a block $A^{i}_{n}=M_{[n,i]}(I_{k_{n,i}})$, we still use $X_{n,i}$ for $Sp(A^{i}_{n})=[0,1]$--that is, $A^{i}_{n}$ is regarded as a homogeneous algebra or a sub-homogeneous algebra over $X_{n,i}$.

\subsection{Definition 2.11} ([EG2], [DG]) Let $X$ be a compact connected space and let $P\in M_{N}(C(X))$ be a projection of rank $n$. The weak variation of a finite set $F\subset PM_{N}(C(X))P$ is defined by
\[
\omega(F)=\sup_{\pi,\sigma}\inf_{u\in U(n)}\max_{a\in F}\|u\pi(a)u^{*}-\sigma(a)\|,
\]
where $\pi,\sigma$ run through the set of irreducible representations of $PM_{N}C(X)P$ into $M_{n}(\mathbb{C})$.

For $F\subset M_{r}(I_{k})$, we define $\omega(F)=\omega(\imath(F))$, where $\imath: M_{r}(I_{k})\longrightarrow M_{rk}(C[0,1])$ is the canonical embedding
and $\imath(F)$ is regarded as a finite subset of $M_{rk}(C[0, 1])$. Let $A$ be a basic $\mathcal{HD}$ block, a finite set $F\subset A$ is said to be weakly approximately constant to within $\varepsilon$ if $\omega(F) <\varepsilon$.

\subsection{2.12} Let $\underline{K}(A)=K_{*}(A)\oplus \bigoplus_{k=2}^{\infty}K_{*}(A,\mathbb{Z}/k)$ be as in [DG]. Let $\Lambda$ be the Bockstein operation between $\underline{K}(A)$s (see 4.1 of [DG]). It is well known that
\[
K_{*}(A,\mathbb{Z}\oplus \mathbb{Z}/k)\cong K_{0}(A\otimes C(T_{II,k}\times S^{1})).
\]
As in [DG], let
\[
K_{*}(A,\mathbb{Z}\oplus \mathbb{Z}/k)^{+}\cong K_{0}(A\otimes C(T_{II,k}\times S^{1}))^{+}.
\]
and let $\underline{K}(A)^{+}$ to be the semigroup of $\underline{K}(A)$ generated by $\{K_{*}(A,\mathbb{Z}\oplus\mathbb{Z}/k\mathbb{Z})^{+}, k=2,3,\cdots\}$.

\subsection{2.13} Let $Hom_{\Lambda}(\underline{K}(A),\underline{K}(B))$ be the homomorphism between $\underline{K}(A)$ and $\underline{K}(B)$ compatible with Bockstein operation $\Lambda$. Associativity of Kasparov $KK$-product gives a map
\[
\Gamma: KK(A,B)\rightarrow Hom_{\Lambda}(\underline{K}(A),\underline{K}(B))
\]

Recall that every element $\alpha\in KK(A,B)$ defines a map $\alpha_{*}\in Hom_{\Lambda}(\underline{K}(A),\underline{K}(B))$. That is, it gives a
sequence of homomorphisms
\[
\alpha^{i}:K_{i}(A)\longrightarrow K_{i}(B) ~~ ~~for~~i=0,1,~~~~and~~~\alpha_{k}^{i}:K_{i}(A,\mathbb{Z}/k)\longrightarrow K_{i}(B,\mathbb{Z}/k) ~~~for~~i=0,1,
\]
which are compatible with the Bockstein operation $\Lambda$.

\subsection{2.14}
Let $\underline{K}(A)_{+}$ be defined in 4.6 of [DG].  Notice that, from [DG, 4.7], $\underline{K}(A)^{+}\subset \underline{K}(A)_{+}$, where for any unital C*-algebra $B$, let
\[
KK(A,B)_{+}=\{\alpha\in KK(A,B),\alpha(\underline{K}(A)_{+})\subseteq \underline{K}(B)_{+}\},  \quad KK(A,B)_{D,+}=\{\alpha\in KK(A,B)_{+},\alpha [1_{A}]\leq [1_{A}]\}.
\]

\subsection{2.15} Let $A=PM_{\bullet}C(T_{II,k})P$, where we use $\bullet$ to denote any possible positive integers. From 5.14 of [G3], we know that an element $\alpha\in KK(A,B)$ is completely determined by $\alpha^{0}: K_{0}(C(T_{II,n}))\longrightarrow K_{0}(B)$ and $\alpha^{1}_{k}: K_{1}(C(T_{II,k}),\mathbb{Z}/k)\longrightarrow K_{1}(B,\mathbb{Z}/k)$,  $k=0,1,2,\cdots.$\

For any positive integer $k\geq 2$, denote
\[
K_{(k)}(A)=K_{0}(A)\oplus K_{1}(A)\oplus K_{0}(A,\mathbb{Z}/k)\oplus K_{1}(A,\mathbb{Z}/k) ~~ and ~~K_{(k)}(A)_{+}=\underline{K}(A)_{+}\cap K_{(k)}(A)
\]

Then for $A=PM_{\bullet}C(T_{II,k})P$, an element $\alpha\in KK(A,B)$ is in $KK(A,B)_{+}$ if and only if $\alpha(K_{(k)}(A)_{+})\subset K_{(k)}(B)_{+}$, $k=0,1$. Note that $K_{(k)}(A)_{+}$ is finitely generated.

\subsection{2.16} Note that, for any C*-algebra $A$, it follows from K\"{u}nneth Theorem that
\[
K_{0}(A\otimes C(W_{k}\times S^{1}))=K_{0}(A)\oplus K_{1}(A)\oplus K_{0}(A,\mathbb{Z}/k)\oplus K_{1}(A,\mathbb{Z}/k)
\]
where $W_{k}=T_{II,k}$.

One can choose finite set $\mathcal{P}\subset M_{\bullet}(A\otimes C(W_{k}\times S^{1}))$ of projections such that
\[
\{[P]\in K_{0}(A)\oplus K_{1}(A)\oplus K_{0}(A,\mathbb{Z}/k)\oplus K_{1}(A,\mathbb{Z}/k)|P\in \mathcal{P}(A)\}\subset \underline{K}(A)
\]
generate $K_{(k)}(A)=K_{0}(A)\oplus K_{1}(A)\oplus K_{0}(A,\mathbb{Z}/k)\oplus K_{1}(A,\mathbb{Z}/k)$, where $W_{k}=T_{II,k}$.

\par ~~If we choose finite set $G(\mathcal{P})\subset A$ large enough and $\delta(\mathcal{P})>0$ small enough, then every $G(\mathcal{P})-\delta(\mathcal{P})$ multiplicative contraction $\phi: A\rightarrow B$ determines a map $\phi_{*}: \mathcal{P}\underline{K}(A)\rightarrow \underline{K}(B)$ which is compatible with the Bockstein operation $\Lambda$ (see [GL]). If $A=PM_{\bullet}C(T_{II,k})P$, then it also defines a $KK$ element $[\phi] \in KK(A,B)$.

\par ~~Let $A=PM_{\bullet}C(T_{II,k})P$, recall that for $G\supseteq G(\mathcal{P})$, $\delta \leq \delta(\mathcal{P})$, a $G-\delta$ multiplicative map $\phi: A\rightarrow B$ is called a quasi $\mathcal{P}\underline{K}$ homomorphism, if there is a homomorphism $\psi: A\rightarrow B$ satisfying $\phi_{*}=\psi_{*}: \mathcal{P}\underline{K}(A)\rightarrow \underline{K}(B)$. If $\phi$ is a quasi $\mathcal{P}\underline{K}$ homomorphism, then $[\phi]\in KK(A,B)_{+}$.

\par ~~Suppose that $A=PM_{\bullet}C(T_{II,k})P$ and suppose that $B=\bigoplus\limits_{i=1}^{n}B^{i}$ is a direct sum of basic $\mathcal{HD}$ bulding blocks. Let $\alpha \in KK(A,B)$, then $\alpha \in KK(A,B)_{+}$ if and only if for each $i\in \{1,2,\cdots,n\}$, $\alpha^{i}\in KK(A,B^{i})$ is either zero or $\alpha^{i}(\mathbf{1}_{C(T_{II,k})})>0$ in $K_{0}(B^{i})$. Recall that in [DG], a $KK$ element $\alpha\in KK(A,B=\bigoplus\limits_{i=1}^{n}B^{i})$ is called m-large if $rank(\alpha^{i}(\mathbf{1}_{A}))\geq m\cdot rank(\mathbf{1}_{A})$, for all $i\in \{1,2,\cdots,n\}$. By 5.5 and 5.6 of [DG] if $\alpha \in KK(A,B)_{+}$ is 6-large and $\alpha([\mathbf{1}_{A}])\leq [\mathbf{1}_{B}]$, then $\alpha$ can be realized by a homomorphism from $A$ to $B$.

\subsection{2.17} Recall that the scale of A, denoted by $\sum A$, is a subset of $K_{0}(A)$ consisting of $[p]$, where $p\in A$ is a projection. As all the C*-algebras $A$ in this paper have cancellation of projections, if $A$ has unit $\mathbf{1}_{A}$, then
\[
\sum A=\big\{x\in K_{0}(A); 0\leq x\leq [\mathbf{1}_{A}]\big\}.
\]
 For two C*-algebras $A,B$, by a "homomorphism" $\alpha$ from $(\underline{K}(A),\underline{K}(A)^{+},\sum A)$ to $(\underline{K}(B),\underline{K}(B)^{+},\sum B)$, it means a system of maps
\[
\alpha_{k}^{i}: K_{i}(A,\mathbb{Z}/k)\rightarrow K_{i}(B,\mathbb{Z}/k); \quad i=0,1,\ \ k=0,2,3,4,\cdots
\]
which are compatible with Bockstein operation and $\alpha=\bigoplus_{k,i}\alpha_{k}^{i}$ satisfies $\alpha(\underline{K}(A)^{+})\subseteq \underline{K}(B)^{+}$ and finally $\alpha_{0}^{0}(\sum A)\subseteq \sum B$.

\subsection{2.18} For a unital C*-algebra $A$, let $TA$ denote the space of tracial states of $A$, i.e. $\tau \in TA$, if and only if $\tau$ is a positive linear map from $A$ to the complex plane $\mathbb{C}$, with $\tau(xy)=\tau(yx)$ and $\tau(\mathbf{1}_{A})=1$. $AffTA$ is the Banach space of all the continuous affine maps from $TA$ to $\mathbb{C}$. (In most references, $AffTA$ is defined to be the set of all the affine maps from $TA$ to $\mathbb{R}$. Our $AffTA$ is a complexification of the standard $AffTA$.) The
element $\mathbf{1}\in AffTA$, defined by $\mathbf{1}(\tau)=1$, for all $\tau \in TA$, is called the unit of $AffTA$. $AffTA$, together with the positive cone $AffTA_{+}$ and the unit element $\mathbf{1}$ forms a
scaled ordered complex Banach space. (Notice that for any element $x\in AffTA$, there are $x_{1}, x_{2}, x_{3}, x_{4}\in AffTA_{+}$ such that $x=x_{1}-x_{2}+ix_{3}-ix_{4}$.)

There is a natural homomorphism $\rho_{A}:K_{0}(A)\rightarrow AffTA$ defined by $\rho_{A}([p])(\tau)=\sum_{i=1}^{n}\tau(p_{ii})$ for $\tau\in TA$ and $[p]\in K_{0}(A)$ represented by projection $p=(p_{i,j})\in M_{n}(A)$.

Any unital homomorphism $\phi: A\rightarrow B$ induces a continuous affine map $T\phi: TB\rightarrow TA$, It turns out that $T\phi$ induces a unital positive linear map
\[
AffT\phi: AffTA\rightarrow AffTB.
\]
If $\phi: A \rightarrow B$ is not unital, we still have the positive linear map
\[
AffT\phi: AffTA\rightarrow AffTB.
\]
but it will not preserve the unit $\mathbf{1}$, only has property $AffT\phi(\mathbf{1}_{AffTA})\leq \mathbf{1}_{AffTB}$.

\subsection{2.19} If $\alpha:(\underline{K}(A),\underline{K}(A)^{+},\sum A)\rightarrow (\underline{K}(B),\underline{K}(B)^{+},\sum B)$ is a homomorphism as in 2.17, then for each projection $p\in A$, there is a projection $q\in B$ such that $\alpha[p] = [q]$.

Notice that for all the C*-algebras $A$ considered in this paper, the following is true: if $p_{1}$, $p_{2}$ are projections and $[p_{1}]=[p_{2}]$ in $K_{0}(A)$, then there is $u\in A$ with $up_{1}u^{*}=p_{2}$. Therefore both $T(pAp)$ and $T(qBq)$ depend only on the class $[p]\in K_{0}(A)$ and $[q]\in K_{0}(B)$. For the classes $[p]\in \sum A(\subset K_{0}(A))$, $T(pAp)$ is taken as part of the invariant $Inv^{0}(A)$. For two classes $[p]\in \sum A$, $[q]\in \sum B$, with $\alpha([p])=[q]$, we will consider the system of continuous affine maps $\xi_{p,q} : T(qBq) \rightarrow T(pAp)$. Such system of maps is said to be compatible if for any two projections $p_{1}\leq p_{2}$ with $\alpha([p_{1}])=[q_{1}]$, $\alpha([p_{2}])=[q_{2}]$ and $q_{1}\leq q_{2}$, the following diagram is commutative:
\[
\CD
  AffT(p_{1}Ap_{1}) @>\xi^{p_{1},q_{1}}>> AffT(q_{1}Bq_{1}) \\
  @V \imath_{T} VV @V \imath_{T} VV \\
  AffT(p_{2}Ap_{2}) @>\xi^{p_{2},q_{2}}>> AffT(q_{2}Bq_{2})
\endCD
\]
where the horizontal maps are induced by $\xi^{p_{1},q_{1}}$ and $\xi^{p_{2},q_{2}}$ respectively, and the
vertical maps are induced by the inclusions $p_{1}Ap_{1}\rightarrow p_{2}Ap_{2}$, $q_{1}Bq_{1}\rightarrow q_{2}Bq_{2}$.

\subsection{2.20} We denote $\big(\underline{K}(A),\underline{K}(A)^{+},\sum A, \{T(pAp)\}_{[p]\in \sum A}\big)$ by $Inv^{0}(A)$. By a ``map'' between the invariants $\big(\underline{K}(A),\underline{K}(A)^{+},\sum A, \{T(pAp)\}_{[p]\in \sum A}\big)$  and $\big(\underline{K}(B),\underline{K}(B)^{+},\sum B, \{T(qBq)\}_{[q]\in \sum B}\big)$, we mean a map
\[
\alpha:\Big(\underline{K}(A),\underline{K}(A)^{+},\sum A\Big)\rightarrow \Big(\underline{K}(B),\underline{K}(B)^{+},\sum B\Big)
\]
as in 2.17 and maps $\xi^{p,q}: T(qBq)\rightarrow T(pAp)$ which are compatible as 2.19. We denote this map by
\[
(\alpha,\xi): \Big(\underline{K}(A),\underline{K}(A)^{+},\sum A,\{T(pAp)\}_{[p]\in \sum A}\Big)\rightarrow \Big(\underline{K}(B),\underline{K}(B)^{+},\sum B,\{T(qBq)\}_{[q]\in \sum B}\Big)
\]
or simply
\[
(\alpha,\xi): Inv^{0}(A)\rightarrow Inv^{0}(B)
\]

\subsection{2.21} Let $A$ be a unital C*-algebra. Let $U(A)$ denote the group of unitaries of $A$ and let $U_{0}(A)$ denote the connected component of $\mathbf{1}_{A}$ in $U(A)$. Let $DU(A)$ and $DU_{0}(A)$ denote the commutator subgroups of $U(A)$ and $U_{0}(A)$, respectively. (Recall that the commutator subgroup of a group $G$ is the subgroup generated by all elements of the form $aba^{-1}b^{-1}$, where $a,b\in G$.) One introduces the following metric $D_{A}$ on $U(A)/\overline{DU(A)}$ (see [NT,\S3]). For $u,v\in U(A)/\overline{DU(A)}$
\[
D_{A}(u,v)=\inf\{\|uv^{*}-c\|:c\in \overline{DU(A)}\}
\]
where, on the right hand side of the equation, we use $u$, $v$ to denote any elements in $U(A)$, which represent the elements $u, v \in U(A)/\overline{DU(A)}$.

Denote the extended commutator group $DU^{+}(A)$, which is generated by $DU(A)\subset U(A)$ and the set $\big\{e^{2\pi i t p}=e^{2\pi i t}p+(\mathbf{1}-p)\in U(A)\big| t \in \mathbb{R}, p\in A \ is \ a \ projection\big\}$. Let $\widetilde{DU(A)}$ denote the closure of $DU^{+}(A)$. That is, $\widetilde{DU(A)}=\overline{DU^{+}(A)}$

\subsection{2.22} Let $A$ be a unital C*-algebra. Let $AffTA$ and $\rho_{A}:K_{0}(A)\rightarrow AffTA$ be as defined as in 2.18.

For simplicity, we will use $\rho K_{0}(A)$ to denote the set $\rho_{A} (K_{0}(A))$. The metric $d_{A}$ on $AffTA/\overline{\rho K_{0}(A)}$ is defined as follows (see [NT, \S3]).

Let $d^{'}$ denote the quotient metric on $AffTA/\overline{\rho K_{0}(A)}$, i.e, for $f,g\in AffTA/\overline{\rho K_{0}(A)}$,
\[
d^{'}(f,g)=\inf\{\|f-g-h\|,\ h\in \overline{\rho K_{0}(A)}\}.
\]
Define $d_{A}$ by
\[
d_{A}(f,g)=
\begin{cases}
2, & \     if \ \ d^{'}(f,g)\geq \frac{1}{2}, \\
|e^{2\pi i d^{'}(f,g)}-1|, & \ if \ \ d^{'}(f,g) < \frac{1}{2}.
\end{cases}
\]
Obviously, $d_{A}(f,g)\leq 2\pi d^{'}(f,g)$.

Let $\widetilde{\rho K_{0}(A)}$ denote the closed real vector space spanned by $\overline{\rho K_{0}(A)}$. That is,
\[
\widetilde{\rho K_{0}(A)}:=\overline{\Big\{\sum\lambda_{i}\phi_{i}\Big|\lambda_{i}\in \mathbb{R}, \phi_{i}\in \rho K_{0}(A) \Big\}}.
\]

For any $u,v\in U(A)/\widetilde{DU(A)}$, define
\[
\overline{D_{A}}(u,v)=\inf\{\|uv^{*}-c\|:c\in \widetilde{DU(A)}\}.
\]
Let $\widetilde{d}^{'}$ denote the quotient metric on $AffTA/\widetilde{\rho K_{0}(A)}$, that is,
\[
\widetilde{d}^{'}(f,g)=\inf\{\|f-g-h\|,\ h\in \widetilde{\rho K_{0}(A)}\}, \quad f,g\in AffTA/\widetilde{\rho K_{0}(A)}.
\]
Define $\widetilde{d}_{A}$ by
\[
\widetilde{d}_{A}(f,g)=
\begin{cases}
2, & \     if \ \ \widetilde{d}^{'}(f,g)\geq \frac{1}{2}, \\
|e^{2\pi i \widetilde{d}^{'}(f,g)}-1|, & \ if \ \ \widetilde{d}^{'}(f,g) < \frac{1}{2}.
\end{cases}
\]

\subsection{2.23} Let
\[
\widetilde{SU(A)}:=\overline{\Big\{x\in U(A)\Big| x^{n}\in \widetilde{DU(A)}\ \ for \ some \ n\in \mathbb{Z}_{+}\setminus\{0\}\Big\}}.
\]
and for simpler, denote $\widetilde{\rho K_{0}}(A)\triangleq \widetilde{\rho K_{0}(A)}$, $\widetilde{DU}(A)=\widetilde{DU(A)}$, $\widetilde{SU}(A)=\widetilde{SU(A)}$.

A unital homomorphism $\phi: A\rightarrow B$ induces a contractive group homomorphism
\[
\phi^{\natural}: U(A)/\widetilde{SU}(A)\longrightarrow U(B)/\widetilde{SU}(B).
\]
If $\phi$ is not unital, then the map $\phi^{\natural}: U(A)/\widetilde{SU}(A)\longrightarrow U(\phi(\mathbf{1}_{A})B\phi(\mathbf{1}_{A}))/\widetilde{SU}(\phi(\mathbf{1}_{A})B\phi(\mathbf{1}_{A}))$ is induced by the corresponding unital homomorphism. In this case, $\phi$ also induces the map $\imath_{*}\circ \phi^{\natural}: U(A)/\widetilde{SU}(A)\longrightarrow U(B)/\widetilde{SU}(B)$, which is denoted by $\phi_{*}$ to avoid confusion.

\subsection{2.24} In [GJL1] and [GJL2], denote
\[
\Big(\underline{K}(A),\underline{K}(A)^{+}, \sum A, \big\{AffT(pAp)\big\}_{[p]\in \sum A}, \big\{U(pAp)/\widetilde{SU}(pAp)\big\}_{[p]\in \sum A } \Big)
\]
by $Inv(A)$.

By a map from $Inv(A)$ to $Inv(B)$, one means
\[
\alpha: \Big(\underline{K}(A),\underline{K}(A)^{+}, \sum A\Big)\longrightarrow \Big(\underline{K}(B),\underline{K}(B)^{+}, \sum B\Big)
\]
as in 2.17, and for each pair $([p],[\overline{p}])\in \sum A\times \sum B$ with $\alpha([p])=[\overline{p}]$, there are an associate unital positive linear map
\[
\xi^{p,\overline{p}}: AffT(pAp)\longrightarrow AffT(\overline{p}B\overline{p})
\]
and an associate contractive group homomorphism
\[
\chi^{p,\overline{p}}: U(pAp)/\widetilde{SU}(pAp)\longrightarrow U(\overline{p}B\overline{p})/\widetilde{SU}(\overline{p}B\overline{p})
\]
satisfying the following compatibility conditions:\\
(a) If $p<q$, then the diagrams
\[
\qquad \qquad \qquad \qquad \qquad \CD
  AffT(pAp) @>\xi^{p,\overline{p}}>> AffT(\overline{p}B\overline{p}) \\
  @V \imath_{T} VV @V \imath_{T} VV       \ \ \qquad \qquad  \qquad \qquad  \qquad            (I) \\
  AffT(qAq) @>\xi^{q,\overline{q}}>> AffT(\overline{q}B\overline{q})
\endCD
\]
and
\[
\qquad \qquad \qquad \quad \quad \CD
  U(pAp)/\widetilde{SU}(pAp) @>\chi^{p,\overline{p}}>> U(\overline{p}A\overline{p})/\widetilde{SU}(\overline{p}B\overline{p}) \\
  @V \imath_{*} VV @V \imath_{*} VV       \ \ \quad \quad  \qquad \quad  \qquad            (II) \\
  U(qAq)/\widetilde{SU}(qAq) @>\chi^{q,\overline{q}}>> U(\overline{q}A\overline{q})/\widetilde{SU}(\overline{q}B\overline{q})
\endCD
\]
commutes, where the vertical maps are induced by inclusions.\\
(b) The following diagram commutes
\[
\qquad \qquad \qquad \qquad \qquad \qquad \ \ \CD
  K_{0}(pAp) @>\rho>> AffT(pAp) \\
  @V \alpha VV @V \xi^{p,\overline{p}} VV       \ \qquad \qquad  \qquad \qquad  \qquad            (III) \\
  K_{0}(\overline{p}B\overline{p}) @>\rho>> AffT(\overline{p}B\overline{p})
\endCD
\]
and therefore $\xi^{p,\overline{p}}$ induces a map (still denoted by $\xi^{p,\overline{p}}$):
\[
\xi^{p,\overline{p}}: AffT(pAp)/\widetilde{\rho K_{0}}(pAp)\longrightarrow AffT(\overline{p}B\overline{p})/\widetilde{\rho K_{0}}(\overline{p}B\overline{p})
\]
(The commutativity of (III) follows from the commutativity of (I), by 1.20 of [Ji-Jiang]. So this is not an
extra requirement.)\\
(c) The following diagrams
\[
\qquad \quad \qquad \quad \quad \CD
  AffT(pAp)/\widetilde{\rho K_{0}}(pAp) @>>> U(pAp)/\widetilde{SU}(pAp) \\
  @V \xi^{p,\overline{p}} VV @V \chi^{p,\overline{p}} VV       \ \ \quad \qquad  \quad \quad  \quad            (IV) \\
  AffT(\overline{p}B\overline{p})/\widetilde{\rho K_{0}}(\overline{p}B\overline{p}) @> >> U(\overline{p}B\overline{p})/\widetilde{SU}(\overline{p}B\overline{p})
\endCD
\]
and
\[
\qquad \qquad \qquad \quad \quad \CD
  U(pAp)/\widetilde{SU}(pAp) @> >> K_{1}(pAp)/torK_{1}(pAp) \\
  @V \chi^{p,\overline{p}} VV @V \alpha_{1} VV       \ \ \quad \quad  \qquad \quad  \quad            (V) \\
  U(\overline{p}B\overline{p})/\widetilde{SU}(\overline{p}B\overline{p}) @>>> K_{1}(\overline{p}B\overline{p})/torK_{1}(\overline{p}B\overline{p})
\endCD
\]
commute, where $\alpha_{1}$ is induced by $\alpha$.

We will denote the map from $Inv(A)$ to $Inv(B)$ by
\[
\begin{split}
(\alpha,\xi,\chi): \Big(\underline{K}(A),\underline{K}(A)^{+}, \sum A, \big\{AffT(pAp)\big\}_{[p]\in \sum A}, \big\{U(pAp)/\widetilde{SU}(pAp)\big\}_{[p]\in \sum A } \Big)  \longrightarrow\\
 \Big(\underline{K}(B),\underline{K}(B)^{+}, \sum B, \big\{AffT(\overline{p}B\overline{p})\big\}_{[\overline{p}]\in \sum B}, \big\{U(\overline{p}B\overline{p})/\widetilde{SU}(\overline{p}B\overline{p})\big\}_{[\overline{p}]\in \sum B } \Big)
 \end{split}
\]

Note that $Inv^{0}$ is part of $Inv(A)$.

\section{Local approximation lemma}

\subsection{Proposition 3.1} Let $A=PM_{\bullet}C(T_{II,k})P$ or $M_{l}(I_{k})$ and $\mathcal{P}$ be as in 2.16. For any finite set $F\subset A$, $\varepsilon >0$, there exist a finite set $G\subset A$ ($G\supset G(\mathcal{P})$ large enough), a positive number $\delta >0$($\delta<\delta(\mathcal{P})$ small enough) such that the following statement is true:

If $B\in \mathcal{HD}$, $\phi,\psi\in Map(A,B)$ are $G-\delta$ multiplicative and $\phi_{*}=\psi_{*}: \mathcal{P}\underline{K}(A)\rightarrow \underline{K}(B)$, then there is a homomorphism $\nu\in Hom(A,M_{L}(B))$ defined by point evaluations and there is a unitary $u\in M_{L+1}(B)$ such that $\|(\phi\oplus\nu)(a)-u(\psi\oplus\nu)(a)u^{*}\|<\varepsilon$ for all $a\in F$

\subsection{Proof} For the special case that $A=PM_{\bullet}C(T_{II,k})P$ and $B$ is homogeneous, this is [G3, Theorem 5.18]. The proof of the general case is completely same (see 2.16 above). Note that, the calculation of $K_{0}(\Pi^{\infty}_{n=1}B_{n}/\bigoplus^{\infty}_{n=1}B_{n})$, $K_{0}(\Pi^{\infty}_{n=1}B_{n}/\bigoplus^{\infty}_{n=1}B_{n},\mathbb{Z}/k)$ and $K_{1}(\Pi^{\infty}_{n=1}B_{n}/\bigoplus^{\infty}_{n=1}B_{n},\mathbb{Z}/k)$ in 5.12 of [G3]
works well for the case that some of $B_{n}$ being of the form $M_{l}(I_{k})$.\\
\qed

\subsection{Lemma 3.2} Let $A=PM_{\bullet}(C(X))P$ or $M_{l}(I_{k})$. Let $F\subset A$ be approximately constant to within $\varepsilon$ (i.e., $\omega(F)<\varepsilon$). Then for any two homomorphisms $\phi,\psi: A\rightarrow B$ defined by point evaluations with $K_{0}\phi=K_{0}\psi$, there exists a unitary $u\in B$ such that $\|\phi(f)-u\psi(f)u^{*}\|<2\varepsilon$ $\forall f\in F$.

\subsection{Proof} The case $A=PM_{\bullet}(C(X))P$ is Lemma 3.5 of [GJLP2]. For the case $A=M_{l}(I_{k})$, one can also find a homomorphism $\phi^{'}: M_{l}(\mathbb{C})\rightarrow B$ such that
\[
\|\phi(f)-\phi^{'}(f(\underline{0}))\|<\varepsilon \qquad \forall f\in F,
\]
where $f(\underline{0})$ as in 2.9. Then one follows the same argument as the proof of Lemma 3.5 of [GJLP2] to get the result.\\
\qed

\subsection{Lemma 3.3} Let $A=PM_{\bullet}(C(X))P$ or $M_{l}(I_{k})$, $\varepsilon >0$, finite set $F\subset A$ with $\omega(A)<\varepsilon$. There exist a finite set $G\subset A$ ($G\supset G(\mathcal{P})$), a number $\delta>0$ ($\delta<\delta(\mathcal{P})$)and a positive integer $L$ such that the following statement is true. If $B\in \mathcal{HD}$, and $\phi,\psi\in Map(A,B)_{1}$ are $G-\delta$ multiplicative and $\phi_{*}=\psi_{*}: \mathcal{P}\underline{K}(A)\rightarrow \underline{K}(B)$ and $\nu: A\rightarrow M_{\infty}(B)$ is a homomorphism defined by point evaluations with $\nu([\mathbf{1}_{A}])\geq L\cdot [\mathbf{1}_{B}]\in K_{0}(B)$, then there is a unitary $u\in (\mathbf{1}_{B}\oplus \nu(\mathbf{1}_{A}))M_{\infty}(B)(\mathbf{1}_{B}\oplus \nu(\mathbf{1}_{A}))$ such that
\[
\|(\phi\oplus\nu)(f)-u(\psi\oplus\nu)(f)u^{*}\| <5\varepsilon \qquad  \forall f\in F.
\]

\subsection{Proof} Let $L_{1}$ be as in Proposition 3.1 and $L=2L_{1}$. Since $\nu([\mathbf{1}_{A}])\geq L\cdot [\mathbf{1}_{B}]=2L_{1}\cdot [\mathbf{1}_{B}]\in K_{0}(B)$, there is a projection $Q< \nu(\mathbf{1}_{A})$ such that $Q$ is equivalent to $\mathbf{1}_{M_{L_{1}}(B)}$. By Proposition 3.1, there exist a homomorphism $\nu_{1}:A\rightarrow QM_{\infty}(B)Q$ defined by point evaluation and a unitary $w\in (\mathbf{1}_{B}\oplus \nu_{1}(\mathbf{1}_{A}))M_{\infty}(B)(\mathbf{1}_{B}\oplus \nu_{1}(\mathbf{1}_{A}))$ such that $\|(\phi\oplus\nu_{1})(f)-w(\psi\oplus\nu_{1})(f)w^{*}\| <\varepsilon \quad  \forall f\in F.$

Since $\nu_{1}$ (and $\nu$ respectively) is homotopic to a homomorphism factoring through $M_{rank(P)}(\mathbb{C})$ (for $A=PM_{\bullet}(C(X))P$) or factoring through $M_{l}(\mathbb{C})$ (for $A=M_{l}(I_{k})$), there is a unital homomorphism $\nu_{2}: A\rightarrow (\nu(\mathbf{1}_{A})-\nu_{1}(\mathbf{1}_{A}))M_{\infty}(B)(\nu(\mathbf{1}_{A})-\nu_{1}(\mathbf{1}_{A}))$ such that $K_{0}(\nu_{1}\oplus\nu_{2})=K_{0}(\nu)$.

By Lemma 3.2, $\nu$ is approximately unitarily equivalent to $\nu_{1}\oplus\nu_{2}$ to within $2\varepsilon$ on $F$. Hence $\phi\oplus \nu$ is approximately unitarily equivalent to $\psi\oplus\nu $ on $F$ to within $2\varepsilon+\varepsilon+2\varepsilon=5\varepsilon$.\\
\qed

\subsection{Lemma 3.4} Let $A=PM_{\bullet}C(T_{II,k})P$, and $\mathcal{P}\subset M_{\bullet}(A\otimes C(W_{k}\times S^{1}))$ as in 2.16. And let $B=\lim\limits_{n\rightarrow \infty}(B_{n},\psi_{nm})$ be a unital simple inductive limit of direct sums of $\mathcal{HD}$ building blocks. Let $\phi:A\rightarrow B$ be a
unital homomorphism. It follows that for any $G\supset G(\mathcal{P})$ and $\delta < \delta(\mathcal{P})$, there is a $B_{n}$ and a unital $G-\delta$ multiplicative contraction $\phi_{1}:A\rightarrow B_{n}$ which is a quasi $\mathcal{P}\underline{K}$-homomorphism such that $\|\psi_{n,\infty}\circ \phi_{1}(g)-\phi(g)\|<\delta$

\subsection{Proof} There is a finite set $G_{1}$ and $\delta_{1}>0$ such that if a complete positive linear map $\psi:A\rightarrow C$ ($C$ is a C*-algebra) and a homomorphism $\phi:A\rightarrow C$ satisfy that $\| \phi(g)-\psi(g)\|<\delta_{1}$ for all $g\in G_{1}$, then $\psi$ is $G-\frac{\delta}{2}$ multiplicative (see Lemma 4.40 in [G3]). Now if $\phi^{'}:A\rightarrow B_{l}$ satisfies $\|\psi_{l,\infty}\circ \phi^{'}(g)-\phi(g)\|<\delta$ for all $g\in G_{1}$, then $\psi_{l,\infty}\circ \phi^{'}$ is $G-\frac{\delta}{2}$ multiplicative. Hence for some $n > l$ (large enough), $\phi_{1}=\psi_{i,n}\circ \phi^{'}:A\rightarrow B_{n}$ is $G-\delta$ multiplicative. By replacing $G$ by $G\cup G_{1}$ and $\delta$ by $min(\delta,\delta_{1})$, we only need to construct a unital quasi $\mathcal{P}\underline{K}$-homomorsim $\phi_{1}:A\rightarrow B_{n}$ such that
\begin{eqnarray*}
 (*)\qquad  \qquad  \qquad  \qquad  \qquad  \qquad  \|\psi_{n,\infty}\circ \phi_{1}(g)-\phi(g)\|<\delta, \qquad  \forall g\in G.  \qquad  \qquad  \qquad  \qquad  \qquad
\end{eqnarray*}\
Furthermore, if $\phi_{1}$ is a $G-\delta$ multiplicative map satisfying the above condition, then
\[
KK(\phi_{1})\times KK(\psi_{n,\infty})=KK(\phi)\in KK(A,B)_{+}.
\]
Note that $K_{(k)}(A)_{+}$ is finitely generated, there is an $m > n$ such that
\[
KK(\phi_{1})\times KK(\psi_{n,m_{1}})\in KK(A,B_{m_{1}})_{+}.
\]
Since $B$ is simple, for certain $m > m_{1}$, $KK(\psi_{n,m}\circ\phi_{1})$ is 6-large and therefore can be realized by a homomorphism. That is, replacing $\phi_{1}$ by $\psi_{n,m}\circ\phi_{1}$, we get a quasi $\mathcal{P}\underline{K}$ -homomorphism. Thus the proof of the lemma is reduced to the construction of $\phi_{1}$ to satisfy $(*)$.

Since $B=\lim\limits_{n\rightarrow \infty}(B_{n},\psi_{nm})$, there is a $B_{n}$ and finite set $F\subset B_{n}$ such that $G \subset_{\frac{\delta}{3}} \psi_{n,\infty}(F)$. Since $B_{n}(\subset B)$ is a nuclear C*-algebra, there are two complete positive contractions $\lambda_{1} : B_{n} \rightarrow M_{N}(\mathbb{C})$ and $\lambda_{2} : M_{N}(\mathbb{C}) \rightarrow B_{n}$ such that
\begin{eqnarray*}
 \|\lambda_{2}\circ \lambda_{1}(f)-f\|<\frac{\delta}{3} \qquad  \forall f\in F.
\end{eqnarray*}
Since $B_{n}$ is a subalgebra of $B$, by Arverson¡¯s Extension Theorem, one can extend the map $\lambda_{1}:B_{n}\rightarrow M_{N}(\mathbb{C})$ to $\beta_{1}:B\rightarrow M_{N}(\mathbb{C})$ such that $\beta_{1}\circ \psi_{n,\infty}=\lambda_{1}$.

One can verify that $\phi_{1}=\lambda_{2}\circ \beta_{1}\circ \phi: A\rightarrow B_{n}$ satisfies the condition $(*)$ as below. For $g \in G$, there is an $f \in F$ such that $\|\phi(g)-\psi_{n,\infty}(f)\|<\frac{\delta}{3}$, and therefore
\begin{eqnarray*}
\|\psi_{n,\infty}\circ \phi_{1}(g)-\phi(g)\|
&   =   &   \|(\psi_{n,\infty}\circ \lambda_{2}\circ \beta_{1} \circ\phi_{1})(g)-\phi(g)\|\\
& \leq  &   \|\psi_{n,\infty}(\lambda_{2}\circ \beta_{1} \circ\phi_{1}(g))-\psi_{n,\infty}(\lambda_{2}\circ \beta_{1} \circ\psi_{n,\infty}(f))\|\\
&   +   & \|(\psi_{n,\infty}\circ \lambda_{2}\circ \lambda_{1} (f)-\psi_{n,\infty}(f)\|+\|(\psi_{n,\infty}(f)-\phi(g)\|  \\
& \leq  &    \frac{\delta}{3}+ \frac{\delta}{3}+ \frac{\delta}{3}=\delta.
\end{eqnarray*}
This ends the proof of the lemma.\\
\qed\

\subsection{Lemma 3.5} For any finite set $F\subset PM_{\bullet}C(T_{II,k})P\triangleq B$ with $\omega(F)< \varepsilon$, there are a finite set $G\supset G(\mathcal{P})$, a positive number $\delta < \delta(\mathcal{P})$
and a positive integer $L$ (this $L$ will be denoted by $L(F, \varepsilon)$ later), such that if $C$ is a $\mathcal{HD}$ basic building block, $p, q \in C$ are two projections with $p+q=\mathbf{1}_{C}$, and $\phi_{0}:B\rightarrow pCp$ and $\varphi_{1}:B\rightarrow qCq$ are two maps satisfying the following conditions:\\
(1) $\phi_{0}$ is quasi $\mathcal{P}\underline{K}$ homomorphism and $G-\delta$ multiplicative, and $\phi_{1}$ is defined by point evaluations (or equivalently, factoring through a finite dimensional C*-algebra), and\\
(2) $rank(q) \geq L ~rank(p)$,\\
then there is a homomorphism $\phi: B \rightarrow C$ such that $\|\phi_{0}\oplus\phi_{1}(f)-\phi(f)\| <5\varepsilon$ for all $f \in F$.

\subsection{Proof} Since $\phi_{0}$ is a quasi $\mathcal{P}\underline{K}$-homomorphism and is $G-\delta$ multiplicative, there is a homomorphism $\phi_{0}^{'}:B\rightarrow pCp$ such that $\phi^{'}_{0*}=\phi_{0*}: \mathcal{P}\underline{K}(B) \rightarrow \underline{K}(C)$. By Lemma 3.3, there is unitary $u \in C$ such that
\begin{eqnarray*}
 \|(\phi_{0}\oplus\phi_{1})(f)-u(\phi_{0}^{'}\oplus\phi_{1})(f)u^{*}\| <5\varepsilon \qquad  \forall f\in F.
\end{eqnarray*}
The homomorphism $\phi=Adu^{*}\circ(\phi_{0}^{'}\oplus\phi_{1})$ is as desired.\\
\qed\\

\subsection{Lemma 3.6} Let $\varepsilon_{1}>\varepsilon_{2}>\cdots>\varepsilon_{n}>\cdots$ be a sequence with $\sum \varepsilon_{i}<+\infty$. Let $A$ be a simple AH algebra with no dimension growth (as in [EGL] and [Li4]). Then $A$ can be written as an $\mathcal{AHD}$ inductive limit
\begin{eqnarray*}
A_{1}=\bigoplus\limits_{i=1}^{t_{1}}A^{i}_{1}\longrightarrow A_{2}=\bigoplus\limits_{i=1}^{t_{1}}A^{i}_{2}\longrightarrow\cdots
\end{eqnarray*}
with the unital subalgebras $B_{n}(=\bigoplus B^{i}_{n})\subset A_{n}(=\bigoplus A^{i}_{n})$ (that is, $B^{i}_{n}\subset A^{i}_{n}$) and $F_{n}\subset B_{n}$ and $G_{n}\subset A_{n}$ with $F_{n}=\bigoplus F^{i}_{n}\subset G_{n}=\bigoplus G^{i}_{n}$ such that the following statements hold:

(1) If $A^{i}_{n}$ is not of type $T_{II}$, then $B^{i}_{n}=A^{i}_{n}$ and $F^{i}_{n}=F^{i}_{n}$. If $A^{i}_{n}$ is of type $T_{II}$, then $B^{i}_{n}=PA^{i}_{n}P\oplus D^{i}_{n}\subset A^{i}_{n}$ with $F^{i}_{n}=\pi_{0}(F^{i}_{n})\oplus \pi_{1}(F^{i}_{n})\subset G^{i}_{n}$ and $\omega(\pi_{0}(F^{i}_{n}))<\varepsilon$, where $D^{i}_{n}$ is a direct sum of the $\mathcal{HD}$ building blocks
other than type $T_{II}$, and $\pi_{0}:B^{i}_{n}\rightarrow PA^{i}_{n}P\triangleq B^{0,i}_{n}$ and $\pi_{1}:B^{i}_{n}\rightarrow D^{i}_{n}$ are canonical projections;

(2) $G^{i}_{n}$ generates $A^{i}_{n}$, $\phi_{n,n+1}(A_{n})\subset B_{n+1}$, $\phi_{n,n+1}(G_{n})\subset F_{n+1}$, and $\overline{\bigcup\limits_{n=1}^{\infty}\phi_{n,\infty}(G_{n})}=\overline{\bigcup\limits_{n=1}^{\infty}\phi_{n,\infty}(F_{n})}=$unit ball of $A$;

(3) Suppose that both $A^{i}_{n}$ and $A^{j}_{n+1}$ are of type $T_{II}$ and $\phi\triangleq \pi_{0}\circ \phi^{i,j}_{n,n+1}|_{B^{0,i}_{n}}:B^{0,i}_{n}\rightarrow B^{0,j}_{n+1}$. Then $\phi(\mathbf{1}_{B^{0,i}_{n}})=p_{0}\oplus p_{1}\in B^{0,j}_{n+1}$ and $\phi=\phi_{0}\oplus\phi_{1}$ with $\phi_{0}\in Hom(B^{0,i}_{n},p_{0}B^{0,j}_{n+1}p_{0})_{1}$, $\phi_{1}\in Hom(B^{0,i}_{n},p_{1}B^{0,j}_{n+1}p_{1})_{1}$ such that $\phi_{1}$ is defined by point evaluations (or equivalently, $\phi_{1}(B^{0,i}_{n})$ is a finite dimensioned sub-algebra of $p_{1}B^{0,j}_{n+1}p_{1}$) and
\begin{eqnarray*}
[p_{1}]\geq L(\pi_{0}(F^{i}_{n}),\varepsilon_{n})\cdot [p_{0}],
\end{eqnarray*}
where $L(\pi_{0}(F^{i}_{n}),\varepsilon_{n})$ is as in Lemma 3.5 (note that $\omega(\pi_{0}(F^{i}_{n}))<\varepsilon_{n}$).

\subsection{Proof} From [Li4], we know that $A$ can be written as an inductive limit of the direct sums of $\mathcal{HD}$ building blocks $A=\lim\limits_{n\rightarrow \infty}(\widetilde{A_{n}},\psi_{nm})$. In our construction, we will choose $A_{n}=\widetilde{A_{k_{n}}}$ and homomorphisms $\phi_{n,n+1}:A_{n}=\widetilde{A_{k_{n}}}\rightarrow A_{n+1}=\widetilde{A_{k_{n+1}}}$ satisfying $KK(\phi_{n,n+1})=KK(\psi_{k_{n},k_{n+1}})$ and $AffT\phi_{n,n+1}$ is close to $AffT\psi_{k_{n},k_{n+1}}$ within any pregiven small number on a pregiven finite set, in such a way that $\phi_{n,n+1}$ also
satisfies the desired condition as in the lemma with certain choices of subalgebras $B_{n}\subset A_{n}$ and and finite subsets $F_{n}\subset B_{n}$ and $G_{n}\subset A_{n}$.

Suppose that we already have $A_{n}=\widetilde{A_{k_{n}}}$ with the unital sub-algebra $B_{n}=\oplus_{i}B^{i}_{n}=\oplus_{i}(B^{0,i}_{n}\oplus D^{i}_{n})\subset A_{n}=\oplus_{i}A^{i}_{n}$, and two subsets $F_{n}\subset B_{n}$, $G_{n}\subset A_{n}$. as required in the lemma. We will construct, roughly, the
next algebra $A_{n+1}=\widetilde{A_{k_{n+1}}}$, the sub-algebra $B_{n+1}\subset A_{n+1}$, two subsets $F_{n+1}\subset B_{n+1}$, $G_{n+1}\subset A_{n+1}$, and the homomorphism $\phi_{n,n+1}:A_{n}\rightarrow A_{n+1}$ to satisfy all the requirements in the lemma.

Applying the decomposition theorem-[G3, Theorem 4.37], as in the proof of the main theorem in [Li4],
for $l>k_{n}$ large enough, and for each block $A_{n}^{i}=\widetilde{A^{i}_{k_{n}}}$ of type $T_{II}$, the homomorphism $\psi^{i,j}_{k_{n},l}$ can be decomposed into three parts $\phi^{1}_{0}\oplus \phi_{1}\oplus\phi_{2}$, roughly described as below: There are mutually orthogonal projections $Q_{0}$,$Q_{1}$,$Q_{2}\in \widetilde{A^{j}_{l}}$ with $\psi^{i,j}_{k_{n},l}(\mathbf{1}_{A^{i}_{n}})=Q_{0}+Q_{1}+Q_{2}$, there are two homomorphisms $\phi_{k}\in Hom(A^{i}_{n},Q_{k}\widetilde{A^{j}_{l}}Q_{k})_{1}$, $(k=1,2)$ and a sufficient multiplicative quasi $\mathcal{P}\underline{K}$ homomorphism $\phi_{k}\in Map(A^{i}_{n},Q_{0}\widetilde{A^{j}_{l}}Q_{0})_{1}$ possessing the following properties:

(i) $\phi^{1}_{0}\oplus \phi_{1}\oplus\phi_{2}$ is close to $\psi^{i,j}_{k_{n},l}$ to within $\varepsilon_{n}$ on $G_{n}^{i}$;

(ii) $\phi^{1}_{0}(\mathbf{1}_{B^{0,i}_{n}})$ is a projection and
\begin{eqnarray*}
[Q_{1}]\geq L(\pi_{0}(F^{i}_{n}),\varepsilon)\cdot [Q_{0}] \quad and \quad \phi_{1}(\mathbf{1}_{B^{0,i}_{n}})\geq L(\pi_{0}(F^{i}_{n}),\varepsilon)\cdot [\phi^{1}_{0}(\mathbf{1}_{B^{0,i}_{n}})];
\end{eqnarray*}

(iii) $\phi_{1}$ is defined by point evaluation at a dense enough finite subset of $Sp(\widetilde{A^{i}_{k_{n}}})$ such that for any homomorphism $\phi_{0}\in Hom(A^{i}_{n},Q_{0}\widetilde{A^{j}_{l}}Q_{0})_{1}$, we have $\omega(\phi_{0}\oplus \phi_{1}(G^{i}_{n}))<\varepsilon_{n+1}$;

(iv) $\phi_{2}$ factors through an interval algebra $D_{1}^{i}$ as $\phi_{2}=\xi_{2}^{i}\circ \xi_{1}^{i}: A^{i}_{n}\stackrel{\xi_{1}^{i}}{\longrightarrow} D_{1}^{i} \stackrel{\xi_{2}^{i}}{\longrightarrow} \widetilde{A^{j}_{l}}$.

Then we can choose $k_{n+1}=l$ and $A_{n+1}=\widetilde{A_{k_{n+1}}}$. Note that $\phi_{0}^{1}$ is a quasi $\mathcal{P}\underline{K}$ homomorphism, one can choose $\phi_{0}: \widetilde{A^{i}_{k_{n}}}\rightarrow \phi_{0}^{1}(\mathbf{1}_{\widetilde{A^{i}_{k_{n}}}})\widetilde{A^{i}_{k_{n+1}}}\phi_{0}^{1}(\mathbf{1}_{\widetilde{A^{i}_{k_{n}}}})$ such that $KK(\phi^{1}_{0}\oplus \phi_{1}\oplus\phi_{2})=KK(\phi_{0}\oplus \phi_{1}\oplus\phi_{2})$ with $\phi^{1}_{0}(\mathbf{1}_{B^{0,i}_{n}})=\phi_{0}(\mathbf{1}_{B^{0,i}_{n}})$. Modify $\psi_{k_{n},k_{n+1}}$ by replacing $\psi^{i,j}_{k_{n},k_{n+1}}$ by $\phi_{0}\oplus \phi_{1}\oplus\phi_{2}$ to define $\phi_{n,n+1}:A_{n}\rightarrow A_{n+1}$.

Fix $j$ with $A^{j}_{n+1}=A^{j}_{k_{n+1}}$ being of type $T_{II}$. Let $I_{0}=\{~i~ |~A^{i}_{n}=\widetilde{A^{j}_{k_{n}}}~is~of~type~T_{II}\}$, $I_{1}=\{~i~ |~A^{i}_{n}=\widetilde{A^{j}_{k_{n}}}~is~not~of~type~T_{II}\}$. Let $P_{i}=(\phi_{0}\oplus\phi_{1})(\mathbf{1}_{A^{i}_{n}})\in A^{j}_{n+1}$, let $P=\bigoplus\limits_{i\in I_{0}}P_{i}\in A_{n+1}^{j}$ and $B^{0,j}_{n+1}\triangleq PA^{j}_{n+1}P$, and let $D^{j}_{n+1}=\bigoplus\limits_{i\in I_{0}}\xi^{i}_{2}(D_{1}^{i})\oplus\bigoplus\limits_{i\in I_{1}}\widetilde{\phi^{i,j}_{k_{n},k_{n+1}}}(\widetilde{A^{j}_{k_{n}}})\subset\widetilde{A^{j}_{k_{n+1}}}$ which is a direct sum of $\mathcal{HD}$ building blocks other than type $T_{II}$. Let $B^{j}_{n+1}=B^{0,j}_{n+1}\oplus D^{j}_{n+1}\subset A^{j}_{n+1}$ if $A^{j}_{n+1}$ is of type $T_{II}$; and $B^{j}_{n+1}=A^{j}_{n+1}$ if $A^{j}_{n+1}$ is not of type $T_{II}$. Choose $G_{n+1}=\oplus G^{j}_{n+1}(\supset \phi_{n,n+1}(G_{n}))$ sufficiently large. If $A^{j}_{n+1}$ is not of type $T_{II}$, let $F^{j}_{n+1}=G^{j}_{n+1}$. If $A^{j}_{n+1}$ is of type $T_{II}$, let $F^{j}_{n+1}=\phi^{-,j}_{n,n+1}(G_{n})$.

Note that $KK(\psi_{k_{n},k_{n+1}})=KK(\phi_{n,n+1})$. Since the rank of $\phi^{1}_{0}(\mathbf{1}_{\widetilde{A^{i}_{k_{n}}}})=\phi_{0}(\mathbf{1}_{\widetilde{A^{i}_{k_{n}}}})$ is much smaller than the ranks of $\phi_{1}(\mathbf{1}_{\widetilde{A^{i,j}_{k_{n}}}})$ and $\phi_{2}(\mathbf{1}_{\widetilde{A^{i,j}_{k_{n}}}})$, $AffT\phi_{n,n+1}$ is very close to $AffT\psi_{k_{n},k_{n+1}}$. Hence $A^{'}=lim(A_{n,n+1},\phi_{n,m})$ has the same Elliott invariant as $A=lim(\widetilde{A_{k_{n}}},\psi_{k_{n},k_{m}})$. Hence $A^{'}\cong A$.\\
\qed\\

\subsection{Corollary 3.7} Let $A$ be a simple $AH$ algebra with no dimension growth. And let $P_{1},P_{2},\cdots,P_{k}\in A$ be a set of mutually orthogonal projections. Then one can write A as inductive limit $A=\lim(A_{n},\phi_{n,m})$ with mutually orthogonal projections $P^{0}_{1},P^{0}_{2},\cdots,P^{0}_{k}\in A_{1}$ such that for each $i$, $\phi_{1,\infty}(P^{0}_{i})=P_{i}$ and
\begin{eqnarray*}
P_{i}AP_{i}=\lim\limits_{n\rightarrow\infty}\Big(\phi_{1,n}(P^{0}_{i})A_{n}\phi_{1,n}(P^{0}_{i}),\phi_{n,m}|_{\phi_{1,n}(P^{0}_{i})A_{n}\phi_{1,n}(P^{0}_{i})}\Big)
\end{eqnarray*}
satisfies the properties of Lemma 3.6.

\subsection{Proof} In the proof of Lemma 3.6, one can assume that there are $P^{'}_{1},P^{'}_{2},\cdots,P^{'}_{k}\in \widetilde{A_{1}}$ with $\psi_{1,\infty}(P^{'}_{i})=P_{i}$ for all $i=1,2,\cdots,k$. Then the construction can be carried out to get our conclusion. Note that, applying
Lemma 1.6.8 of [G3], one can strengthen [G3, Theorem 4.37] such that the following is true: For a set of
pre-given orthogonal projections $p_{1},p_{2},\cdots,p_{k}\in A_{n}$, one can further require that $\psi_{0}\in Map(A_{n},Q_{0}A_{m}Q_{0})$ satisfies that $\psi_{0}(p_{i})$ are projections and $(\psi_{0}\oplus\psi_{1}\oplus\psi_{2})(p_{i})=\phi_{n,m}(p_{i})$ for $i=1,2,\cdots,k$.\\
\qed\\

The following lemma is the main technique lemma of this section.

\subsection{Lemma 3.8} Let $A$,$A^{'}$ be simple $\mathcal{AHD}$ inductive limit algebras with $A_{1}\stackrel{\phi_{1,2}}{\longrightarrow} A_{2} \stackrel{\phi_{2,3}}{\longrightarrow} \cdots \stackrel{}{\longrightarrow} A_{n}\cdots \stackrel{}{\longrightarrow} A$ and $A^{'}_{1}\stackrel{\psi_{1,2}}{\longrightarrow} A^{'}_{2} \stackrel{\psi_{2,3}}{\longrightarrow} \cdots \stackrel{}{\longrightarrow} A^{'}_{n}\cdots \stackrel{}{\longrightarrow} A^{'}$ being described in Lemma 3.6. Let $\Lambda:A\rightarrow A^{'}$ be a homomorphism. Let $F\subset A_{m}$ be a finite set and $\varepsilon >0$. Then there is an $A^{'}_{l}$ and a homomorphism $\Lambda_{1}:A_{m}\rightarrow A^{'}_{l}$ such that
\begin{eqnarray*}
\|\Lambda\circ\phi_{m,\infty}(f)-\psi_{l,\infty}\circ\Lambda_{1}(f)\|<\varepsilon \qquad \forall f\in F.
\end{eqnarray*}

\subsection{Proof} One can choose $n > m$ large enough such that $\phi_{m,n}(F)\subset_{\frac{\varepsilon}{4}}F_{n}$ and $5\varepsilon_{n}<\frac{\varepsilon}{8}$. Note that $\phi_{m,n}(A_{m})\subset B_{n}$, where $B_{n}$ is as in Lemma 3.6. We will construct a homomorphism $\phi:B_{n}\rightarrow A^{'}_{l}$ such that
\begin{eqnarray*}
 (*)\qquad  \qquad  \qquad  \qquad   \|\Lambda\circ\phi_{n,\infty}(f)-\psi_{l,\infty}(\phi(f))\|<\frac{\varepsilon}{2} \qquad \forall f\in F_{n}\subset B_{n}\subset A_{n}. \qquad  \qquad  \qquad
\end{eqnarray*}
Then the homomorphism $\phi\circ\phi_{m,n}$ is as desired.

Let $I_{0}=\big\{~i~ |~B^{i}_{n}~is~of~type~T_{II}\big\}$ and $I_{1}=\big\{~i~|~B^{i}_{n}~is~not~of~type~T_{II}\big\}$. Then $\big\{\mathbf{1}_{B^{0,i}_{n}}\big\}_{i\in I_{0}}\bigcup\big\{\mathbf{1}_{D^{i}_{n}}\big\}_{i\in I_{0}}\bigcup\big\{\mathbf{1}_{B^{i}_{n}}\big\}_{i\in I_{1}}$ are mutually orthogonal projections. Hence $\big\{\Lambda\circ\phi_{n,\infty}(\mathbf{1}_{B^{0,i}_{n}})\big\}_{i\in I_{0}}\bigcup\big\{\Lambda\circ\phi_{n,\infty}(\mathbf{1}_{D^{i}_{n}})\big\}_{i\in I_{0}}\bigcup\big\{\Lambda\circ\phi_{n,\infty}(\mathbf{1}_{B^{i}_{n}})\big\}_{i\in I_{1}}$ are mutually orthogonal projections in $A^{¡ä}$. One can choose $n_{1}$ (large enough) and mutually orthogonal projections $\big\{P_{i}\big\}_{i\in I_{0}}\bigcup\big\{Q_{i}\big\}_{i\in I_{0}}\bigcup\big\{R_{i}\big\}_{i\in I_{1}}\subset A^{'}_{n_{1}}$ and a unitary $u\in A^{'}$ such that $\|u-1\|<\frac{\varepsilon}{16}$ and $\psi_{n_{1},\infty}(P_{i})=u^{*}(\Lambda\circ\phi_{n,\infty}(\mathbf{1}_{B^{0,i}_{n}}))u$,  ~ $\psi_{n_{1},\infty}(Q_{i})=u^{*}(\Lambda\circ\phi_{n,\infty}(\mathbf{1}_{D^{i}_{n}}))u$ for $i\in I_{0}$ and $\psi_{n,\infty}(R_{i})=u^{*}(\Lambda\circ\phi_{n,\infty}(\mathbf{1}_{B^{i}_{n}}))u$ for $i\in I_{1}$.

Note that $\|Adu-id\|<\frac{\varepsilon}{8}$. Replacing $\Lambda: A\rightarrow A^{'}$ by $\Lambda^{'}=Adu\circ\Lambda: A\stackrel{\Lambda}{\longrightarrow}A^{'}\stackrel{Adu}{\longrightarrow}A^{'}$ to make $(*)$ true, it suffices to construct $\phi: B_{n}\rightarrow A^{'}_{l}$ for certain $l>n_{1}$ such that
\begin{eqnarray*}
 (**)\qquad  \qquad  \qquad  \quad   \|\Lambda^{'}\circ\phi_{n,\infty}(f)-\psi_{l,\infty}\circ\phi(f)\|<\frac{3\varepsilon}{8} \qquad \forall f\in F_{n}\subset B_{n}\subset A_{n}. \qquad  \qquad  \quad
\end{eqnarray*}

Such construction can be carried out for each of the blocks $B^{0,i}_{n}$ and $D^{i}_{n}$ for $i\in I_{0}$, $B^{i}_{n}$ for $i\in I_{1}$. Namely, we need to construct
\begin{eqnarray*}
\phi|_{B^{0,i}_{n}}: B^{0,i}_{n}\longrightarrow \psi_{n_{1},l}(P_{i})A^{'}_{l}\psi_{n_{1},l}(P_{i}) \qquad \forall i\in I_{0};
\end{eqnarray*}
\begin{eqnarray*}
\phi|_{D^{i}_{n}}: D^{i}_{n}\longrightarrow \psi_{n_{1},l}(Q_{i})A^{'}_{l}\psi_{n_{1},l}(Q_{i}) \qquad \forall i\in I_{0};
\end{eqnarray*}
\begin{eqnarray*}
\phi|_{B^{i}_{n}}: B^{i}_{n}\longrightarrow \psi_{n_{1},l}(R_{i})A^{'}_{l}\psi_{n_{1},l}(R_{i}) \qquad \forall i\in I_{1}.
\end{eqnarray*}
separately, to satisfy the condition
\begin{eqnarray*}
  \|\Lambda^{'}\circ\phi_{n,\infty}(f)-\psi_{l,\infty}\circ\phi(f)\|<\frac{3\varepsilon}{8}
\end{eqnarray*}
for all $f\in \big\{\pi_{0}(F^{i}_{n})\big\}_{i\in I_{0}}\cup\big\{\pi_{1}(F^{i}_{n})\big\}_{i\in I_{0}}\cup\big\{F^{i}_{n}\big\}_{i\in I_{1}}$.

For the blocks $\big\{D^{i}_{n}\big\}_{i\in I_{0}}$ and $\big\{B^{i}_{n}\big\}_{i\in I_{1}}$, the existence of such homomorphisms follows from the fact that the domain algebras are stably generated. So we only need to construct $\phi: B^{0,i}_{n}\longrightarrow \psi_{n_{1},l}(P_{i})A^{'}_{l}\psi_{n_{1},l}(P_{i})$ for $l$ large enough.

Let $J_{0}=\big\{~j~ |~B^{j}_{n+1}~is~of~type~T_{II}\big\}$ and $J_{1}=\big\{~j~|~B^{i}_{n+1}~is~not~of~type~T_{II}\big\}$. Let $\widetilde{P}^{i,j}=\pi_{0}(\phi^{i,j}_{n,n+1}(\mathbf{1}_{B^{0,i}_{n}}))\in B^{0,j}_{n+1}$ for $j\in J_{0}$, $\widetilde{Q}^{i,j}=\pi_{1}(\phi^{i,j}_{n,n+1}(\mathbf{1}_{B^{0,i}_{n}}))\in D^{j}_{n+1}$ for $j\in J_{0}$, $\widetilde{R}^{i,j}=\pi_{1}(\phi^{i,j}_{n,n+1}(\mathbf{1}_{B^{0,i}_{n}}))\in D^{j}_{n+1}$ for $j\in J_{1}$. (Here we only consider the case $i\in I_{0}$).

As in Lemma 3.6, we have the decomposition $\pi_{0}\circ \phi^{i,j}_{n,n+1}(\mathbf{1}_{B^{0,i}_{n}})=p_{0}+p_{1}\in B^{0,j}_{n+1} $ and $\pi_{0}\circ \phi^{i,j}_{n,n+1}|_{B^{0,i}_{n}}=\phi_{0}\oplus \phi_{1}$ with $\phi_{0}\in Hom(B^{0,i}_{n},p_{0}B^{0,j}_{n+1}p_{0})_{1}$ and $\phi_{1}\in Hom(B^{0,i}_{n},p_{1}B^{0,j}_{n+1}p_{1})_{1}$. Denote $p_{0},p_{1}$ by $p^{i,j}_{0}$ and $p^{i,j}_{1}$, then $\widetilde{P}^{i,j}=p^{i,j}_{0}\oplus p^{i,j}_{1}$. It follows that
\begin{eqnarray*}
\Big\{\Lambda^{'}\big(\phi_{n+1,\infty}(p^{i,j}_{0})\big)\Big\}_{j\in J_{0}}\bigcup \Big\{\Lambda^{'}\big(\phi_{n+1,\infty}(p^{i,j}_{1})\big)\Big\}_{j\in J_{0}}\bigcup \Big\{\Lambda^{'}\big(\phi_{n+1,\infty}(\widetilde{Q}^{i,j})\big)\Big\}_{j\in J_{0}}\bigcup\Big \{\Lambda^{'}\big(\phi_{n+1,\infty}(\widetilde{R}^{i,j})\big)\Big\}_{j\in J_{1}}
\end{eqnarray*}
is a set of mutually orthogonal projections with sum to be $\Lambda^{'}\big(\phi_{n,\infty}(\mathbf{1}_{B^{0,i}_{n}})\big)=\psi_{n_{1},\infty}(P_{i})\in A^{'}$. For $n_{2}>n_{1}$ (large enough), there are mutually orthogonal projections
\begin{eqnarray*}
\{P^{i,j}_{0}\}_{j\in J_{0}}\cup \{P^{i,j}_{1})\}_{j\in J_{0}}\cup \{Q^{i,j}\}_{j\in J_{0}}\cup \{{R}^{i,j}\}_{j\in J_{1}}\subset \psi_{n_{1},n_{2}}(P_{i})A^{'}_{n_{2}}\psi_{n_{1},n_{2}}(P_{i})
\end{eqnarray*}
and a unitary $v\in \psi_{n,\infty}(P_{i})A^{'}\psi_{n,\infty}(P_{i})$ such that $\|v-1\|<\frac{\varepsilon}{16}$, and
\begin{eqnarray*}
\psi_{n_{2},\infty}(P^{i,j}_{0})=v^{*}\Big(\Lambda^{'}\big(\phi_{n+1,\infty}(p^{i,j}_{0})\big)\Big)v, \\
\psi_{n_{2},\infty}(P^{i,j}_{1})=v^{*}\Big(\Lambda^{'}\big(\phi_{n+1,\infty}(p^{i,j}_{1})\big)\Big)v, \\
\psi_{n_{2},\infty}(Q^{i,j})=v^{*}\Big(\Lambda^{'}\big(\phi_{n+1,\infty}(\widetilde{Q}^{i,j})\big)\Big)v
\end{eqnarray*}
for all $j\in J_{0}$; and
\begin{eqnarray*}
\psi_{n_{2},\infty}(R^{i,j})=v^{*}\Big(\Lambda^{'}\big(\phi_{n+1,\infty}(\widetilde{R}^{i,j})\big)\Big)v
\end{eqnarray*}
for all $j\in J_{1}$.

Let $\widetilde{\Lambda}=Adv\circ \Lambda^{'}$. Then the construction of $\phi|_{B^{0,i}_{n}}$ satisfying $(**)$ for all $f \in \pi_{0}(F_{0}^{i})$ is reduced to the construction of homomorphisms
\begin{eqnarray*}
\xi^{j}_{0}:  B^{0,j}_{n}\longrightarrow \psi_{n_{2},l}(P^{i,j}_{0}\oplus P^{i,j}_{1})A^{'}_{l}\psi_{n_{2},l}(P^{i,j}_{0}\oplus P^{i,j}_{1}),\\
\xi^{j}_{1}:  \widetilde{Q}^{i,j}D^{j}_{n+1}\widetilde{Q}^{i,j}\longrightarrow \psi_{n_{2},l}(Q^{i,j})A^{'}_{l}\psi_{n_{2},l}(Q^{i,j}) \quad \
\end{eqnarray*}
for all $j\in J_{0}$; and
\begin{eqnarray*}
\xi^{j}:  \widetilde{R}^{i,j}B^{j}_{n+1}\widetilde{R}^{i,j}\longrightarrow \psi_{n_{2},l}(R^{i,j})A^{'}_{l}\psi_{n_{2},l}(R^{i,j})
\end{eqnarray*}
for all $j\in J_{1}$, such that
\begin{eqnarray*}
(***)\quad  \quad      \|\widetilde{\Lambda}\circ\phi_{n,\infty}(f)-\psi_{l,\infty}\circ\xi^{j}_{0}(f)\|<\frac{\varepsilon}{4} \quad \forall f \in \pi_{0}(F_{n}^{i})\subset B_{n}^{0,i}\quad and \ j\in J_{0}. \qquad \qquad\qquad   \qquad\\
  \  \  \    \|\widetilde{\Lambda}\circ\phi_{n+1,\infty}(f)-\psi_{l,\infty}\circ\xi^{j}_{1}(f)\|<\frac{\varepsilon}{4} \quad \forall f \in \pi_{1}\Big(\phi^{i,j}_{n,n+1}|_{B^{0,i}_{n}}(\pi_{0}(F_{n}^{i}))\Big)\subset \widetilde{Q}^{i,j}D^{j}_{n+1}\widetilde{Q}^{i,j} \quad and\ j\in J_{0}. \quad  \quad\\
  \quad  \quad  \quad    \|\widetilde{\Lambda}\circ\phi_{n+1,\infty}(f)-\psi_{l,\infty}\circ\xi^{j}(f)\|<\frac{\varepsilon}{4} \quad \forall f \in \phi^{i,j}_{n,n+1}|_{B^{0,i}_{n}}(\pi_{0}(F_{n}^{i}))\subset \widetilde{R}^{i,j}B^{j}_{n+1}\widetilde{R}^{i,j}\quad and \ j\in J_{1}.  \ \quad  \quad
\end{eqnarray*}
(Warning: The domain of $\xi_{0}^{j}$ is $B^{0,i}_{n}$ which is a sub-algebra of $B_{n}$ but not a sub-algebra of $B_{n+1}$. On the other hand, $\xi^{j}_{1}(j\in J_{0})$ and $\xi^{j}(j\in J_{1})$ are homomorphisms from $\widetilde{Q}^{i,j}D^{j}_{n+1}\widetilde{Q}^{i,j}(j\in J_{0})$ and $\widetilde{R}^{i,j}B^{j}_{n+1}\widetilde{R}^{i,j}(j\in J_{1})$ which are subalgebras of $B_{n+1}$).

The existence of the homomorphisms $\xi^{j}_{1}(j\in J_{0})$ and $\xi^{j}(j\in J_{1})$ follows from the fact that the corresponding domain algebras $\widetilde{Q}^{i,j}D^{j}_{n+1}\widetilde{Q}^{i,j}(j\in J_{0})$ and $\widetilde{R}^{i,j}B^{j}_{n+1}\widetilde{R}^{i,j}(j\in J_{1})$ are stably generated¡ªof course, we need to choose $l>n_{2}$ large enough.

So we only need to construct $\xi^{j}_{0}$ to satisfy $(***)$ above. Let $G \supset G(\mathcal{P})$ and $\delta< \delta(\mathcal{P})$ be as in Lemma 3.5 for $\pi_{0}(F_{n}^{i})$ and $\varepsilon_{n}$ (note that $\omega(\pi_{0}(F_{0}^{i}))<\varepsilon_{n}$). Recall that $L(\pi_{0}(F^{i}_{n}), \varepsilon_{n})$ is also from Lemma 3.5. Recall that $\phi_{0}\in Hom(B^{0,i}_{n},p^{i,j}_{0}B^{j}_{n+1}p^{i,j}_{0})$, $\phi_{1}\in Hom(B^{0,i}_{n},p^{i,j}_{1}B^{j}_{n+1}p^{i,j}_{1})$, and that $\phi_{1}(B^{0,i}_{n})$ is a finite dimensional algebra. There is a homomorphism $\lambda_{1}:\phi_{1}(B^{0,i}_{n})\rightarrow \psi_{n_{2},l}(P^{i,j}_{1})A^{'}_{l}\psi_{n_{2},l}(P^{i,j}_{1})$ (for $l$ large enough) such that
\begin{eqnarray*}
\|\psi_{l,\infty}\circ\lambda_{1}(f)-\widetilde{\Lambda}(\phi_{n+1,\infty}(f))\|<\frac{\varepsilon}{16} \qquad \forall f \in \phi_{1}(\pi_{0}(F_{n}^{i})).
\end{eqnarray*}

Applying Lemma 3.4 to the inductive limit, $\lim\big(\psi_{n_{2},m}(P^{i,j}_{0})A^{'}_{m}\psi_{n_{2},m}(P^{i,j}_{0}),\psi_{m,l}\big)$, the finite set $\phi_{0}(\pi_{0}(F_{n}^{i}))\subset p^{i,j}_{0}B^{0,j}_{n+1}p^{i,j}_{0}$ and the homomorphism $\widetilde{\Lambda}\circ \phi_{n+1,\infty}:p^{i,j}_{0}B^{0,j}_{n+1}p^{i,j}_{0}\rightarrow \psi_{n_{2},\infty}(P^{i,j}_{0})A^{'}_{m}\psi_{n_{2},\infty}(P^{i,j}_{0})$, for $l(>n_{2})$ large enough, one can obtain a $\phi_{0}(G)-\delta$ multiplicative quasi-$\mathcal{P}\underline{K}$ homomorphism $\lambda_{0}:p^{i,j}_{0}B^{0,j}_{n+1}p^{i,j}_{0}\rightarrow \psi_{n_{2},l}(P^{i,j}_{0})A^{'}_{m}\psi_{n_{2},l}(P^{i,j}_{0})$ such that
\begin{eqnarray*}
\|\psi_{l,\infty}\circ\lambda_{0}(f)-(\widetilde{\Lambda}\circ\phi_{n+1,\infty})(f)\|<\frac{\varepsilon}{16} \qquad \forall f \in \phi_{0}(\pi_{0}(F_{n}^{i})).
\end{eqnarray*}

Let $\xi^{'}=\lambda_{0}\circ\phi_{0}\oplus\lambda_{1}\circ\phi_{1}:B^{0,i}_{n}\longrightarrow \psi_{n_{2},l}(P^{i,j}_{0}\oplus P^{i,j}_{1})A^{'}_{l}\psi_{n_{2},l}(P^{i,j}_{0}\oplus P^{i,j}_{1})$. Then
\begin{eqnarray*}
(****)\qquad  \qquad \qquad \quad    \|\widetilde{\Lambda}\circ\phi_{n,\infty}(f)-\psi_{l,\infty}\circ\xi^{'}(f)\|<\frac{\varepsilon}{8} \qquad \forall f \in \pi_{0}(F_{n}^{i}). \qquad \qquad \qquad  \qquad
\end{eqnarray*}

On the other hand, since $[P^{i,j}_{1}]\geq L(\pi_{0}(F^{i}_{n}),\varepsilon_{n})\cdot [P^{i,j}_{0}]$ in K-theory, and $\lambda_{1}\circ\phi_{1}$ is a homomorphism with finite dimensional image, we know that $\xi^{'}=\lambda_{0}\circ\phi_{0}\oplus\lambda_{1}\circ\phi_{1}$ satisfies the condition of Lemma 3.6--note that $\lambda_{0}$ is $\phi_{0}(G)-\delta$ multiplicative implies that $\lambda_{0}\circ\phi_{0}$ is $G-\delta$ multiplicative. By Lemma 3.6, there is a homomorphism $\xi^{j}_{0}:B^{0,j}_{n}\longrightarrow \psi_{n_{2},l}(P^{i,j}_{0}\oplus P^{i,j}_{1})A^{'}_{l}\psi_{n_{2},l}(P^{i,j}_{0}\oplus P^{i,j}_{1})$ such that $\|\xi_{0}^{j}(f)-\xi^{'}(f)\|\leq 5\varepsilon_{n}<\frac{\varepsilon}{8}$. Combining with $(****)$ we know that $\xi^{j}_{0}$ satisfies $(***)$ as desired, and therefore the lemma is proved.\\
\qed

\subsection{Definition 3.9} A C*-algebra $A$ is said to have the ideal property if each closed two-sided ideal in $A$ is generated by its projections.

\subsection{Remark 3.10} All simple, unital C*-algebras have the ideal property. The direct sum of TAI algebras have the ideal property. The inductive limit of C*-algebras with the ideal property have the ideal property.

\par~~

The following results is due to Gong, Jiang, Li.

\subsection{Proposition 3.11} ([GJL1]) Suppose that $A=\lim(A_{n},\phi_{n,m})$ and $B=\lim(B_{n},\psi_{n,m})$ are two (not necessarily unital) $\mathcal{AHD}$ inductive limit algebras with the ideal property. Suppose that there is an isomorphism
\[
\alpha: \Big(\underline{K}(A),\underline{K}(A)^{+}, \sum A\Big)\longrightarrow \Big(\underline{K}(B),\underline{K}(B)^{+}, \sum B\Big)
\]
which is compatible with Bockstein operations. Suppose that for each projection $p \in A$ and $\overline{p} \in B$ with $\alpha([p])=[\overline{p}]$, there exist a unital positive linear isomorphism
\[
\xi^{p,\overline{p}}: AffT(pAp)\longrightarrow AffT(\overline{p}B\overline{p})
\]
and an isometric group isomorphism
\[
\gamma^{p,\overline{p}}: U(pAp)/\widetilde{SU}(pAp)\longrightarrow U(\overline{p}B\overline{p})/\widetilde{SU}(\overline{p}B\overline{p})
\]
satisfying the following compatibility conditions:\\
(1) For each pair of projections $p<q<\in A$ and $\overline{p}<\overline{q}<\in B$ with $\alpha([p])=[\overline{p}]$, $\alpha([q])=[\overline{q}]$, the diagrams
\[
\CD
  AffT(pAp) @>\xi^{p,\overline{p}}>> AffT(\overline{p}B\overline{p}) \\
  @V  VV @V  VV  \\
  AffT(qAq) @>\xi^{q,\overline{q}}>> AffT(\overline{q}B\overline{q})
\endCD
\]
and
\[
\CD
  U(pAp)/\widetilde{SU}(pAp) @>\gamma^{p,\overline{p}}>> U(\overline{p}A\overline{p})/\widetilde{SU}(\overline{p}B\overline{p}) \\
  @V VV @V  VV \\
  U(qAq)/\widetilde{SU}(qAq) @>\gamma^{q,\overline{q}}>> U(\overline{q}A\overline{q})/\widetilde{SU}(\overline{q}B\overline{q})
\endCD
\]
commute, where the vertical maps are induced by the inclusion homomorphisms.\\
(2) The maps $\alpha_{0}$ and $\xi^{p,\overline{p}}$ are compatible, that is, the diagram
\[
\CD
  K_{0}(pAp) @>\rho>> AffT(pAp) \\
  @V \alpha_{0} VV @V \xi^{p,\overline{p}} VV \\
  K_{0}(\overline{p}B\overline{p}) @>\rho>> AffT(\overline{p}B\overline{p})
\endCD
\]
commutes (this is not an extra requirement, since it follows from the commutativity of the first diagram in
(1) above by [Ji-Jiang]), and then we have the map (still denoted by $\xi^{p,\overline{p}}$):
\[
\xi^{p,\overline{p}}: AffT(pAp)/\widetilde{\rho K_{0}}(pAp)\longrightarrow AffT(\overline{p}B\overline{p})/\widetilde{\rho K_{0}}(\overline{p}B\overline{p})
\]\\
(3) The map $\xi^{p,\overline{p}}$ and $\gamma^{p,\overline{p}}$ are compatible, that is, the diagram
\[
\CD
  AffT(pAp)/\widetilde{\rho K_{0}}(pAp) @> \widetilde{\lambda}_{A}>> U(pAp)/\widetilde{SU}(pAp) \\
  @V \xi^{p,\overline{p}} VV @V \gamma^{p,\overline{p}} VV \\
  AffT(\overline{p}B\overline{p})/\widetilde{\rho K_{0}}(\overline{p}B\overline{p}) @>\widetilde{\lambda}_{B}>> U(\overline{p}B\overline{p})/\widetilde{SU}(\overline{p}B\overline{p})
\endCD
\]
commutes.\\
(4) The map $\alpha_{1}: K_{1}(pAp)/torK_{1}(pAp)\longrightarrow K_{1}(\overline{p}B\overline{p})/torK_{1}(\overline{p}B\overline{p})$ (note that $\alpha$ keeps the positive cone of $\underline{K}(A)^{+}$ and therefore takes $K_{1}(pAp)\subset K_{1}(A)$ to $K_{1}(\overline{p}B\overline{p})\subset K_{1}(B)$ is compatible with $\gamma^{p,\overline{p}}$, that is, the diagram
\[
\CD
  U(pAp)/\widetilde{SU}(pAp) @>\widetilde{\pi}_{pAp}>> K_{1}(pAp)/torK_{1}(pAp) \\
  @V \gamma^{p,\overline{p}} VV @V \alpha_{1} VV \\
  U(\overline{p}B\overline{p})/\widetilde{SU}(\overline{p}B\overline{p}) @>\widetilde{\pi}_{\overline{p}B\overline{p}}>> K_{1}(\overline{p}B\overline{p})/torK_{1}(\overline{p}B\overline{p})
\endCD
\]
commutes.\\
Then there is an isomorphism $\Gamma: A\rightarrow B$ such that\\
(a)$K(\Gamma)=\alpha$, and \\
(b)If $\Gamma_{p}: pAp\rightarrow \Gamma(p)B\Gamma(p)$ is the restriction of $\Gamma$, then $AffT(\Gamma_{p})=\xi^{p,\overline{p}}$ and $\Gamma^{\natural}_{p}=\gamma^{p,\overline{p}}$, where $[\overline{p}]=[\Gamma(p)]$.

\subsection{Proposition 3.10} ([GJL2, Proposition 2.38]) Let $A,B \in \mathcal{HD}$ or $\mathcal{AHD}$ be unital C*-algebras. Suppose that $K_{1}(A)=tor(K_{1}(A))$ and $K_{1}(B)=tor(K_{1}(B))$. It follows that $Inv^{0}(A)\cong Inv^{0}(B)$ implies that $Inv(A)\cong Inv(B)$.

\section{Main Theorems}

\subsection{Theorem 4.1} If $A$ is an ATAI C*-algebra, then $A$ is an $\mathcal{AHD}$ algebra with ideal property.

\subsection{Proof} By remark 3.10, $A$ has the ideal property. So we will only prove $A$ is an $\mathcal{AHD}$ algebra.

Suppose that $A$ is the inductive limit $(A_{n}=\bigoplus\limits_{i=1}^{t_{n}}A_{n}^{i},\phi_{n,m})$, where $A^{i}_{n}$ are simple $\mathcal{AHD}$ algebras. Since all $A^{i}_{n}$ are simple, without lose of generality, we can assume that all the homorphisms $\phi_{n,m}$ are injective. We will construct a sequence of sub-C*-algebras $B^{i}_{n}\subset A^{i}_{n}$ which are direct sums of $\mathcal{HD}$ building blocks and homomorphisms $\psi_{n,n+1}: B_{n}=\bigoplus B_{n}^{i}\rightarrow B_{n+1}=\bigoplus B_{n+1}^{i}$ such that the diagram

$$
\xymatrix{
B_{1} \ar@{^(->}[dd]^{\imath_{1}} \ar[rr]^{\psi_{1,2}} & & B_{2}  \ar@{^(->}[dd]^{\imath_{2}} \ar[rr]  & & \cdots \\
& & & &\\
A_{1} \ar[rr]^{\phi_{1,2}}  & & A_{2} \ar[rr] & & \cdots
}
$$

is approximately commutative in the sense of Elliott and $\overline{\bigcup\limits_{n=1}^{\infty}\phi_{n,\infty}(\imath_{n}(B_{n}))}=A$.

Let $\{a_{n,k}\}^{\infty}_{k=1}$ be a dense subset of the unit ball of $A_{n}$ and $\varepsilon_{n}=\frac{1}{2^{n}}$. Let $G_{1}\subset A_{1}$ be defined by $G_{1}=\bigoplus\limits_{i=1}^{t_{1}}\{\pi_{i}(a_{11})\}\subset\bigoplus\limits_{i=1}^{t_{1}}A_{1}^{i}=A_{1}$ where $\pi_{1}:A_{1}\rightarrow A^{i}_{1}$ are canonical projections.

Fix $i\in\{1,2,\cdots,t_{1}\}$, one can write $A^{i}_{1}=\lim\limits_{n\rightarrow\infty}(C_{n}, \xi_{n,m})$ as in Lemma 3.6 with injective homomorphisms $\xi_{n,m}$.

For $\pi_{i}(G_{1})\subset A^{i}_{1}$, we choose $n$ large enough such that $\pi_{i}(G_{1})\subset_{\varepsilon_{1}} \xi_{n,\infty}(C_{n})$. And let $B^{i}_{1}=C_{n}$ and $F^{i}_{1}$ a finite set with $\pi_{i}(G_{1})\subset_{\varepsilon_{1}} \imath(F_{1}^{i})$, where $\imath$ is the inclusion homomorphism $\xi_{n,\infty}$. Let $B=\bigoplus\limits_{i=1}^{t_{1}}B_{1}^{i}\hookrightarrow \bigoplus\limits_{i=1}^{t_{1}}A_{1}^{i}$, and $F_{1}=\bigoplus F_{1}^{i}$.

We will construct the following diagram

$$
\xymatrix@R=0.1ex{
F_{1}=\bigoplus F_{1}^{i} \ & &  F_{2}=\bigoplus F_{2}^{i} \ \ \ \   & &F_{n}=\bigoplus F_{n}^{i} \  \\
\bigcap  \ & & \bigcap \ \ & &\bigcap}
$$
\vspace{-3mm}
$$
\xymatrix{
& B_{1}=\bigoplus B_{1}^{i} \ar@{^(->}@<1mm>[d]_{\imath_{1}} \ar[rr]^{\psi_{1,2}} & & B_{2}=\bigoplus B_{2}^{i} \ar@{^(->}@<1mm>[d]_{\imath_{2}} \ar[rr] & & \cdots B_{n}=\bigoplus B_{n}^{i} \ar@{^(->}@<3mm>[d]_{\imath_{n}} \ar[r] &\\
& A_{1}=\bigoplus A_{1}^{i} \ar[rr]^{\phi_{1,2}} & & A_{2}=\bigoplus A_{2}^{i} \ar[rr] & & \cdots A_{n}=\bigoplus A^{i}_{n} \ar[r] &
}
$$
\vspace{-3mm}
$$
\xymatrix@R=0.1ex{
 \bigcup &  & \bigcup \ \ \ \ &  &  \bigcup  \ \ \\
\ \ \ G_{1}=\bigoplus G_{1}^{i} \  & &  G_{2}=\bigoplus G_{2}^{i} \ \ \ \ & & G_{n}=\bigoplus G_{n}^{i}
}
$$

such that

(1) $G_{k}\subset_{\varepsilon_{k}}\imath_{k}(F_{k})$, $G_{k}\supset \phi_{k-1,k}(G_{k-1}\cup\{a_{k-1,k}\})\cup\{a_{k,j}\}^{k}_{j=1}$;

(2) $\|\imath_{k}\circ\psi_{k-1,k}(f)-\phi_{k-1,k}\circ\imath_{k-1}(f)\|<\varepsilon_{k-1}, \qquad \forall f\in F_{k-1}$.

The construction will be carried out by induction. Suppose that we have the diagram until $F_{n}\subset B_{n}\hookrightarrow A_{n}\supset G_{n}$. We will construct the next piece of the diagram.

Let $P^{i,j}=\phi^{i,j}_{n,n+1}(\mathbf{1}_{A^{i}_{n}})\in A^{j}_{n+1}$. Then $\{P^{i,j}\}^{t_{n}}_{i=1}$ is a set of mutually orthogonal projections in $A^{j}_{n+1}$. Apply Corollary 3.7 for $A^{j}_{n+1}$ in place of $A$ and $\{P^{i,j}\}^{t_{n}}_{i=1}$ in places of $\{P_{1},P_{2},\cdots,P_{k}\}$, we can write $A^{j}_{n+1}=\lim(C_{k},\lambda_{k,l})$ with $P^{i,j}A^{j}_{n+1}P^{i,j}=\lim\limits_{k}(Q^{i,j}_{k}C_{k}Q^{i,j}_{k},\lambda_{k,l}|_{Q^{i,j}_{k}C_{k}Q^{i,j}_{k}})$ (where $\lambda_{k,l}(Q^{i,j}_{k})=Q^{i,j}_{l}\in C_{l}$) being $\mathcal{AHD}$ inductive limit algebra as described in Lemma 3.6. (Warning: Do not confuse these $C_{k}$ with the $C_{n}$ in the expression $A^{i}_{1}=\lim\limits_{n\rightarrow \infty}(C_{n},\xi_{n,m})$.)

For each pair $i,j$, we apply Lemma 3.8 to the homomorphism $\phi^{i,j}_{n,n+1}:A^{i}_{n}\rightarrow P^{i,j}A^{j}_{n+1}P^{i,j}$ in place of $\Lambda: A\rightarrow A^{'}$, and $F^{i}_{n}\subset B^{i}_{n}$ in place of $F \in A_{m}$, there is an $l$ (large enough) and a homomorphism $\psi^{i,j}_{n,n+1}:B_{n}^{i}\rightarrow Q^{i,j}_{l}C_{l}Q^{i,j}_{l}$ such that
\begin{eqnarray*}
\|\phi^{i,j}_{n,n+1}\circ\imath_{n}(f)-\imath_{n+1}\circ\psi^{i,j}_{n,n+1}(f)\|<\varepsilon \qquad \forall f \in F_{n}^{i},
\end{eqnarray*}
where $\imath_{n+1}|_{Q^{i,j}_{l}C_{l}Q^{i,j}_{l}}=\lambda_{l,\infty}|_{Q^{i,j}_{l}C_{l}Q^{i,j}_{l}}$.

Finally choose $G^{j}_{n+1}\supset \pi_{j}\big[(G_{n}\cup\{a_{n,n+1}\})\cup\{a_{n+1,j}\}^{n+1}_{j=1}\cup\imath_{n+1}(\bigcup\limits_{k}\psi^{i,j}_{n,n+1}(F^{i}_{n}))\big]$ and $G_{n+1}=\bigoplus G^{j}_{n+1}$. By increasing $l$, we can assume that there is an $F^{'j}_{n+1}\subset C_{l}$ such that $G_{n+1}^{j}\subset_{\varepsilon_{n+1}}F^{'j}_{n+1}$. Define $B^{j}_{n+1}$ to be $C_{l}$ which is a subalgebra of $A^{j}_{n+1}$. Let $F^{j}_{n+1}\subset B^{j}_{n+1}$ be defined by $F^{j}_{n+1}=F^{'j}_{n+1}\cup\big(\bigoplus\limits_{i}\phi^{i,j}_{n,n+1}(F^{i}_{n})\big)$ and let $F_{n+1}=\bigoplus F^{j}_{n+1}$. This ends our inductive construction. By the Elliott intertwining argument, $\lim(B_{n},\psi_{n,m})=\lim(A_{n},\phi_{n,m})$. This ends the proof.\\
\qed

\subsection{Theorem 4.2} Let $A=\lim(A_{n},\phi_{n,m})$ and $B=\lim(B_{n},\psi_{n,m})$ be inductive limit algebras with $A_{n}^{i}$, $B^{i}_{n}$ being unital separable nuclear simple TAI algebras with UCT. Suppose that there is an isomorphism
\[
\alpha: \Big(\underline{K}(A),\underline{K}(A)^{+}, \sum A\Big)\longrightarrow \Big(\underline{K}(B),\underline{K}(B)^{+}, \sum B\Big)
\]
which is compatible with Bockstein operations. Suppose that for each projection $p \in A$ and $\overline{p} \in B$ with $\alpha([p])=[\overline{p}]$, there exist a unital positive linear isomorphism
\[
\xi^{p,\overline{p}}: AffT(pAp)\longrightarrow AffT(\overline{p}B\overline{p})
\]
and an isometric group isomorphism
\[
\gamma^{p,\overline{p}}: U(pAp)/\widetilde{SU}(pAp)\longrightarrow U(\overline{p}B\overline{p})/\widetilde{SU}(\overline{p}B\overline{p})
\]
satisfying the following compatibility conditions:\\
(1) For each pair of projections $p<q<\in A$ and $\overline{p}<\overline{q}<\in B$ with $\alpha([p])=[\overline{p}]$, $\alpha([q])=[\overline{q}]$, the diagrams
\[
\CD
  AffT(pAp) @>\xi^{p,\overline{p}}>> AffT(\overline{p}B\overline{p}) \\
  @V  VV @V  VV  \\
  AffT(qAq) @>\xi^{q,\overline{q}}>> AffT(\overline{q}B\overline{q})
\endCD
\]
and
\[
\CD
  U(pAp)/\widetilde{SU}(pAp) @>\gamma^{p,\overline{p}}>> U(\overline{p}A\overline{p})/\widetilde{SU}(\overline{p}B\overline{p}) \\
  @V VV @V  VV \\
  U(qAq)/\widetilde{SU}(qAq) @>\gamma^{q,\overline{q}}>> U(\overline{q}A\overline{q})/\widetilde{SU}(\overline{q}B\overline{q})
\endCD
\]
commute, where the vertical maps are induced by the inclusion homomorphisms.\\
(2) The maps $\alpha_{0}$ and $\xi^{p,\overline{p}}$ are compatible, that is, the diagram
\[
\CD
  K_{0}(pAp) @>\rho>> AffT(pAp) \\
  @V \alpha_{0} VV @V \xi^{p,\overline{p}} VV \\
  K_{0}(\overline{p}B\overline{p}) @>\rho>> AffT(\overline{p}B\overline{p})
\endCD
\]
commutes (this is not an extra requirement, since it follows from the commutativity of the first diagram in
(1) above by [Ji-Jiang]), and then we have the map (still denoted by $\xi^{p,\overline{p}}$):
\[
\xi^{p,\overline{p}}: AffT(pAp)/\widetilde{\rho K_{0}}(pAp)\longrightarrow AffT(\overline{p}B\overline{p})/\widetilde{\rho K_{0}}(\overline{p}B\overline{p})
\]\\
(3) The map $\xi^{p,\overline{p}}$ and $\gamma^{p,\overline{p}}$ are compatible, that is, the diagram
\[
\CD
  AffT(pAp)/\widetilde{\rho K_{0}}(pAp) @> \widetilde{\lambda}_{A}>> U(pAp)/\widetilde{SU}(pAp) \\
  @V \xi^{p,\overline{p}} VV @V \gamma^{p,\overline{p}} VV \\
  AffT(\overline{p}B\overline{p})/\widetilde{\rho K_{0}}(\overline{p}B\overline{p}) @>\widetilde{\lambda}_{B}>> U(\overline{p}B\overline{p})/\widetilde{SU}(\overline{p}B\overline{p})
\endCD
\]
commutes.\\
(4) The map $\alpha_{1}: K_{1}(pAp)/torK_{1}(pAp)\longrightarrow K_{1}(\overline{p}B\overline{p})/torK_{1}(\overline{p}B\overline{p})$ (note that $\alpha$ keeps the positive cone of $\underline{K}(A)^{+}$ and therefore takes $K_{1}(pAp)\subset K_{1}(A)$ to $K_{1}(\overline{p}B\overline{p})\subset K_{1}(B)$ is compatible with $\gamma^{p,\overline{p}}$, that is, the diagram
\[
\CD
  U(pAp)/\widetilde{SU}(pAp) @>\widetilde{\pi}_{pAp}>> K_{1}(pAp)/torK_{1}(pAp) \\
  @V \gamma^{p,\overline{p}} VV @V \alpha_{1} VV \\
  U(\overline{p}B\overline{p})/\widetilde{SU}(\overline{p}B\overline{p}) @>\widetilde{\pi}_{\overline{p}B\overline{p}}>> K_{1}(\overline{p}B\overline{p})/torK_{1}(\overline{p}B\overline{p})
\endCD
\]
commutes.\\
Then there is an isomorphism $\Gamma: A\rightarrow B$ such that\\
(a)$K(\Gamma)=\alpha$, and \\
(b)If $\Gamma_{p}: pAp\rightarrow \Gamma(p)B\Gamma(p)$ is the restriction of $\Gamma$, then $AffT(\Gamma_{p})=\xi^{p,\overline{p}}$ and $\Gamma^{\natural}_{p}=\gamma^{p,\overline{p}}$, where $[\overline{p}]=[\Gamma(p)]$.

\subsection{Proof} It follows from the fact of Theorem 4.1 that all ATAI algebras are $\mathcal{AHD}$ algebras with ideal property and Proposition 3.11 ([GJL1]) that all $\mathcal{AHD}$ algebras with ideal property are classified by $Inv(A)$.\\
\qed

\subsection{Corollary 4.3} (A)([Jiang1], Theorem 3.11) Let $A=\lim_{n\rightarrow \infty}(A_{n}=\bigoplus A_{n}^{i}, \phi_{n,m})$ and $B=\lim_{n\rightarrow \infty}(B_{n}=\bigoplus B_{n}^{i}, \psi_{n,m})$ be inductive limits with $A_{n}^{i}$, $B_{n}^{i}$ being unital separable nuclear simple TAI algebras with UCT and torsion $K_{1}$-group. Assume that there is an isomorphism $\alpha\in Hom_{\Lambda}(\underline{K}(A),\underline{K}(B))$ such that
\[
\alpha(\underline{K}(A)^{+})=\underline{K}(B)^{+}, \qquad \ \  \alpha(\sum A)=\sum B.
\]
and for each pair of projections $p\in A$, $q\in B$ with $\alpha([p])=[q]$, there is a continuous
affine homomorphism
\[
\xi^{p,q}: T(qBq)\rightarrow T(pAp)
\]
which is compatible in the sense of 2.19, then there is an isomorphism $\Lambda: A \rightarrow B$ which induces $\alpha$ and $\xi$ above.

(B) ([Fa]) Let $A$ and $B$ be two ATAF algebras. Suppose that is an isomorphism of ordered groups
\[
\alpha: \Big(\underline{K}(A),\underline{K}(A)^{+},\sum A \Big) \rightarrow \Big(\underline{K}(B),\underline{K}(B)^{+},\sum B \Big)
\]
which preserves the action of the Bockstein operations. Then there is a $*$--isomorphism $\varphi: A\rightarrow B$ with $\varphi_{*}=\alpha$.

\subsection{Proof} (A) It follows from the fact that $K_{1}(A)=tor K_{1}(A)$, $K_{1}(B)=tor K_{1}(B)$ and  Proposition 3.10 ([GJL2, Proposition 2.38]).

(B) From Theorem 4.1, if $A$ is an ATAF C*-algebra, then $A$ is an $\mathcal{AHD}$ algebra of real rank zero which is classified in [DG] (Note that $\mathcal{AHD}$ algebra in [DG] is denoted by ASH algebra).\\
\qed

\end{document}